\documentclass[12pt]{amsart}
\newcounter{defcounter}
\usepackage{amssymb,amsmath,amscd,graphicx, enumitem,
latexsym,amsthm}
\usepackage{amssymb,latexsym,amsmath,amscd,graphicx, amsxtra,verbatim}
  \usepackage[all]{xy}
  \usepackage{pdfsync}
  \usepackage[mathscr]{euscript}

\theoremstyle{plain}
\newtheorem{theorem}{Theorem}

\newtheorem{proposition}[theorem]{Proposition}

\newtheorem{lemma}[theorem]{Lemma}

\newtheorem{proposition.definition}[theorem]{Proposition/Definition}

\newtheorem{italicclaim}[theorem]{Claim}

\newtheorem{theoremalpha}{Theorem}

\newtheorem{conjecturealpha}[theoremalpha]{Conjecture}

\theoremstyle{definition}
\newtheorem{definition}[theorem]{Definition}

\newtheorem{remark}[theorem]{Remark}
\newtheorem{example}[theorem]{Example}
\newtheorem{claim}[theorem]{Claim}
\newtheorem{question}[theorem]{Question}

\newtheorem{problem}[theorem]{Problem}

 \theoremstyle{remark}

\newcommand{\lra}{\longrightarrow}

\newcommand{\ra}{\rightarrow}
\newcommand{\noi}{\noindent}
\newcommand{\PP}{\mathbf{P}}

\newcommand{\CC}{\mathbf{C}}
\newcommand{\QQ}{\mathbf{Q}}

\newcommand{\OO}{\mathcal{O}}

\newcommand{\CCC}{\mathcal{C}}
\newcommand{\DDD}{\mathcal{D}}
\newcommand{\EEE}{\mathcal{E}}
\newcommand{\XXX}{\mathcal{X}}
\newcommand{\YYY}{\mathcal{Y}}
\newcommand{\ZZZ}{\mathcal{Z}}

\newcommand{\EE}{\mathcal{E}}

\newcommand{\GGG}[2]{\Gamma \big( {#1},{#2}\big)}

\newcommand{\HH}[3]{H^{{#1}} \big( {#2} , {#3}
\big) }

\newcommand{\pr}{\prime}

\newcommand{\dra}{\dashrightarrow}

\newcommand{\Jac}{\textnormal{Jac}}

\newcommand{\Linser}[1]{| \mspace{1.5mu} {#1}
\mspace{1.5mu} |}
\newcommand{\linser}[1]{\Linser{  {#1}  }}

\newcommand{\gon}{\textnormal{gon}}
\newcommand{\pro}{\textnormal{pr}}

\newcommand{\ol}[1]{\overline{#1}}

\newcommand{\GG}{\mathbf{G}}

\newcommand{\Pic}{\textnormal{Pic}}
\newcommand{\prim}{\textnormal{prim}}

\newcommand{\corrdeg}{\textnormal{corr.\,deg}}
\newcommand{\irr}{\textnormal{irr}}
 
\newcommand{\Hom}{\textnormal{Hom}}
\newcommand{\van}{\textnormal{van}}
\newcommand{\Adj}{\textnormal{Adj}}
\newcommand{\HS}{\textnormal{HS}}

\newcommand{\covgon}{\textnormal{cov.\,gon}}
\newcommand{\BVA}[1]{\textnormal{(BVA)}_{#1}}
\newcommand{\covgenus}{\textnormal{cov.\,genus}}

\numberwithin{theorem}{section}
\begin{document}

\title{Measures of Association Between Algebraic Varieties}

 \author{Robert Lazarsfeld}
  \address{Department of Mathematics, Stony Brook University, Stony Brook, New York 11794}
 \email{{\tt robert.lazarsfeld@stonybrook.edu}}
 \thanks{Research  of the first author partially supported by NSF grant DMS-1739285.}

 \author{Olivier Martin}
  \address{Department of Mathematics, Stony Brook University, Stony Brook, New York 11794}
 \email{{\tt olivier.martin@stonybrook.edu}}

\maketitle

\tableofcontents

\setlength{\parskip}{.175in minus .08in}

\section*{Introduction}

The purpose of this paper is to begin the exploration of what we call -- for lack of a better term -- measures of association between two irreducible projective varieties of the same dimension. Roughly speaking,  the idea is to study from various points of view the minimal complexity of  correspondences between them. Here we set up some foundations,  extend  the main results of \cite{BDELU} to correspondences,  and consider joint covering invariants of curves and hypersurfaces.   In addition we propose a number of questions and conjectures that we hope may be of interest. 

To start with the simplest invariant, fix irreducible  complex projective varieties $X$ and $Y$ of dimension $n$. Recall that a correspondence between $X$ and $Y$ is given by an irreducible $n$-fold  \[Z \, \subseteq \, X \times Y
\]
that dominates both factors. One can view the product of the degrees of $Z$ over $X$ and $Y$ as a first measure of its non-triviality.  We  define the \textit{correspondence degree} between $X$ and $Y$ to be minimum of this product over all correspondences:
\[ \corrdeg(X, Y) \ = \ \min_{Z\subseteq X \times Y}
\Big\{  \deg(Z/X) \cdot \deg(Z/Y) \Big \}.\footnote{Here and below it  is to be understood that the minimum is taken over subvarieties $Z \subseteq X \times Y$ as above.}
\]
Thus
 \[
 \corrdeg(X,Y) \, = \, 1 \ \Longleftrightarrow \ X \, \sim_{\text{birat}} \, Y. \] (In fact, $\log( \corrdeg)$ defines a metric on the set of birational equivalence classes of varieties of fixed dimension.)
We wish to study this integer for natural pairs of varieties $(X, Y)$.

 Note to begin with that there is an evident upper bound on $\corrdeg(X, Y)$ in terms of the degrees of irrationality of $X$ and $Y$.\footnote{By definition  the degree of irrationalty $\irr(V)$ of an irreducible projective variety $V$ of dimension $n$ is the least degree of a rational covering $V \dra \PP^n$. } Specifically, by taking the fibre product of rational coverings $X \dra \PP^n$ and $Y \dra \PP^n$ one sees that 
 \[
 \corrdeg(X,Y) \ \le \ \irr(X) \cdot \irr(Y). \tag{*}
 \]
Our intuition is that equality holding in (*) points in the direction of  $X$ and $Y$ being  as ``birationally independent" as possible.  

Suppose for instance that $X$ and $Y$ are very general smooth curves of genera $g_X, g_Y \ge 1$.  It is well known that then  $\Pic(X \times Y) = \Pic(X) \times \Pic(Y)$. This implies immediately (Example \ref{Very.General.Curve.Example}) that
\[   \corrdeg(X,Y) \ = \ \gon(X) \cdot \gon(Y) ,   \]
and that moreover  the minimal correspondences arise as  fibre products of pencils $X \to \PP^1 \ , \ Y \to \PP^1$ computing the gonalities.  

Our first  result is an analogous independence statement for very general hypersurfaces of large degree.
\begin{theoremalpha}   \label{TheoremA}
Let 
$X, Y \subseteq \PP^{n+1}$ 
be very general hypersurfaces of degrees $d , e \ge 2n + 2$. Then
\[
\corrdeg(X,Y) \ = \ \irr(X)\cdot \irr(Y).
\]
Moreover the minimal degree correspondences arise birationally as fibre products of  coverings of projective space.
\end{theoremalpha}
\noi We remark that according to \cite{BDELU}, one has $\irr(X) = d-1$, $\irr(Y) = e-1$. The proof of the Theorem proceeds along the lines of  \cite{BDELU} and \cite{Yang}, except that a geometric argument in those papers  involving  gonality is replaced by a computation with the cohomology of hypersurfaces. We also show (Theorem \ref{Corr.Deg.General.X.Y.Bound}) that if $A$ is a very ample line bundle on an arbitrary smooth projective variety $M$ of dimension $n+1$, and if $X_d \in \linser{dA}$ and  $Y_e \in \linser{eA}$ are very general divisors, then $\corrdeg(X,Y)$ grows like $de$ for $d, e \gg 0$.

Turning to a different line of thought, it has proven fruitful to study the geometry of curves covering a fixed variety. Specifically, given an irreducible projective variety $V$ of dimension $n$, define:
\begin{align*}
\covgenus(V) \ &= \ \min \Bigg  \{ g \ge 0 \ 
\Bigg | \parbox{2.7in}{\begin{center} Given a general point $x \in V$, $\exists$ an irreducible curve $C \subseteq V$ through $x$ with  $p_g(C) = g$. \end{center}}  
\Bigg  \}. \\ \\
\end{align*}
\begin{align*}
\covgon(V) \ &= \ \min \Bigg  \{ c > 0 \ 
\Bigg | \parbox{2.7in}{\begin{center} Given a general point $x \in V$, $\exists$ an irreducible curve $C \subseteq V$ through $x$ with  $\gon(C) = c$. \end{center}}
\Bigg  \}.
\end{align*}
Viewed as a measure of irrationality,  the covering gonality in particular has been the focus of considerable recent activity (eg \cite{Bast}, \cite{BDELU}, \cite{BCFS}).  For a pair of varieties, it seems natural to consider the genus or gonality of families of curves that simultaneously cover both. Specifically, given irreducible $n$-folds  $X, Y$ as above, define
\begin{align*}
\covgenus(X,Y) \ &= \ \min_{ Z\subseteq X \times Y} \big\{ \covgenus(Z) \big \} \\ \covgon(X,Y) \ &= \ \min_{Z\subseteq X \times Y} \big\{ \covgon(Z) \big \}.
\end{align*}
We refer to these respectively as the joint covering genus and gonality of $X$ and $Y$. Evidently
\begin{align*}
\covgenus(X,Y) \ &\ge \ \max \big \{ \, \covgenus(X) \, , \, \covgenus(Y) \, \big \} \\ \covgon(X,Y) \ &\ge \ \max \big \{ \, \covgon(X) \, , \, \covgon(Y) \, \big \}
\end{align*}
and equality means in effect  that there is a family of curves computing the invariant for one of the varieties that  also covers the other.\footnote{Strictly speaking, one should keep in mind the possibility that curves covering $Z$ might not map birationally to curves covering $X$ or $Y$.} In the other direction we prove:\begin{align*}
  \covgon(X, Y) \ &\le \ \covgon(X) \cdot \covgon(Y) \\
  \covgenus(X,Y) \ &\lesssim \ 3 \cdot \covgenus(X) \cdot \covgenus(Y).
\end{align*}

These  invariants    seem  to be non-trivial already when $X$ and $Y$ are smooth curves: in this case we are looking for the minimal genus or gonality of (the smooth model of) a curve $Z$ that covers both $X$ and $Y$. For the joint covering genus we establish
\begin{theoremalpha} \label{Curve.Genus.Bound} Assume that $X$ and $Y$ are very general  curves of genera $g_X, g_Y \ge 1$. Then
\[
g_X g_Y + \textnormal{(linear in $g_X, g_Y$)}  \ \le \ \covgenus(X,Y) \ \le \ \tfrac{5}{4}\cdot  g_X g_Y+ \textnormal{(linear in $g_X, g_Y$)}.
\]
\end{theoremalpha}
\noi As for the covering gonality, it  is elementary that  $\covgon(X,Y)\lesssim (g_Xg_Y)/4$. We propose \begin{conjecturealpha} 
There exists a number $a > 0$ such that if $X, Y$ are very general curves as in Theorem \ref{Curve.Genus.Bound}, then
\[
 \covgon(X, Y) \ \ge \ a \cdot g_X   g_Y   .
\]
\end{conjecturealpha}
\noi Unfortunately we 
 are so far unable to establish any non-trivial lower bounds. In fact, it is only with some difficulty that we    prove:
\begin{theoremalpha} \label{Hyperelliptic.Curve.Proposition}
Let $X$ and $Y$ be very general hyperelliptic curves of genera     $ g_X  \ge 2$ and $g_Y \ge 3$. Then
\[
\covgon(X,Y) \ = \ 4.\]
\end{theoremalpha}
 \noi When $g_X = g_Y = 2$, the joint covering gonality is either $3$ or $4$. We also compute the covering gonality when one or both of the curves is elliptic. 

In higher dimensions, we adapt  arguments of Ein \cite{Ein} and Voisin \cite{Voisin} to derive essentially additive lower bounds for very general hypersurfaces:
\begin{theoremalpha} \label{Hypersurface.Covering.Bound.Intro.Statement}
Let $X, Y \subseteq \PP^{n+1}$ be very general hypersurfaces of degrees $d$ and $e$. Then
\begin{align*}
   \covgon(X, Y) \ &\ge \  d + e - (3n + 2).   \\ 
   \covgenus(X, Y) \  &\ge \ \frac{d^2 + e^2 - 3n^2}{2} \, + \, \textnormal{(lower order terms)}. 
   \end{align*}

\end{theoremalpha}
\noi We conjecture that the joint covering gonality actually grows multiplicatively in $d$ and $e$, and that the genus comes closer to the elementary upper bound
 \begin{equation} \label{Equation.in.Intro}
  \covgenus(X, Y) \ \lesssim \ \frac{de(d+e)}{2} 
  \end{equation}
 indicated in Remark \ref{Bounds.Are.Additive.Remark}.
 
There are several  questions and problems that come to mind about such measures of association. In addition to  those already mentioned, we discuss a few in the final section of the paper. 

 Concerning organization: we start in $\S 1$ with some general remarks about the invariants we consider. The second section is devoted to the proof of Theorem \ref{TheoremA}. In \S 3, we adapt the approach of Ein and Voisin to study joint covering invariants, proving in particular Theorems \ref{Hypersurface.Covering.Bound.Intro.Statement} and \ref{Curve.Genus.Bound}. Degeneration arguments pioneered by Pirola, Alzati and Voisin are used to establish Theorem \ref{Hyperelliptic.Curve.Proposition} in \S 4. Finally, \S 5 is devoted to some open problems and conjectures.

We work throughout over the complex numbers $\CC$. For a   smooth projective variety $V$ we write $K_V$ for a canonical divisor, but  sometimes denote by $\omega_V$ the canonical bundle of $V$, so that $\omega_V = \OO_V(K_V)$. As in equation \eqref{Equation.in.Intro}
 above,  we allow ourselves to write informally  $\lesssim$  or $ \gtrsim $ to indicate that the quantity on the left is bounded by an expression whose dominant terms appear on the right; it will generally be clear how one could arrive at a precise statement. When we say that $X$ and $Y$ are very general curves or hypersurfaces, we mean that the pair $(X,Y)$ is a very general point in the product of the relevant parameter spaces.

We thank Francesco Bastianelli, Nathan Chen, Ciro Ciliberto, Hannah Larson, John Sheridan, David Stapleton and Ruijie Yang for valuable conversations and discussions. We are particularly grateful to  Dima Dudko and Dennis Sullivan for indirectly launching this project. Following a colloquium talk by the first author on measures of irrationality, Dudko suggested that one might view the results of \cite{BDELU} and related papers as measuring the ``distance" between a variety and $\PP^n$, and   he asked whether the  story could then be generalized to study arbitrary pairs of varieties. It was this question, and a resulting discussion with Sullivan, that started us thinking about the matters appearing here.

 \numberwithin{equation}{section}

\section{Measures of association}

This section is devoted to some preliminary remarks about the invariants we will consider.

Let $X$ and $Y$ be smooth complex projective varieties of dimension $n$. Modifying slightly the terminology appearing in the Introduction, by a \textit{correspondence}  between $X$ and $Y$ we will mean a smooth irreducible projective variety $Z$ of dimension $n$ together with a morphism 
\[   u : Z \lra X \times Y  \]
that is birational onto its image. We denote by $a$ and $b$  the projections  from $Z$ to $X$ and $Y$, and we assume that they are dominant:
\begin{equation} \label{Correspondence.Setup}
\begin{gathered}
\xymatrix{
& Z \ar[dl]_a\ar[dr]^b  \\
X & & Y.
 }
 \end{gathered}
\end{equation} 
Note  that $a$ and $b$,  being dominant, are generically finite. The condition that $Z \lra X \times Y$ be birational to its image guarantees that a general fibre of $a$ (or $b$) is naturally realized as a subset of $Y$   (or $X$). Observe finally that a correspondence $Z$ between $X$ and $Y$ is well-defined in the birational category: given alternative birational models $X^\pr$ and $Y^\pr$ of $X$ and $Y$, one can find $Z^\pr \sim_{\text{bir}} Z$ sitting in the analogue of \eqref{Correspondence.Setup}.

\subsection*{Correspondence degree}   The general goal of the present paper is to study the ``minimal complexity" of a correspondence between fixed varieties. A first invariant in this direction involves the degrees of a correspondence:
\begin{definition}
Let $X$ and $Y$ be smooth projective $n$-folds. The \textit{correspondence degree} between $X$ and $Y$ is defined to be
\[
\corrdeg(X, Y) \ = \ \min_Z \big\{ \deg(a) \cdot \deg(b) \big \},
\]
the minimum being taken over all correspondences as in   \eqref{Correspondence.Setup}.
\end{definition}
\vskip -10pt
\noi Thus
 $\corrdeg(X,Y) = 1$ if and only if $X  \sim_{\text{bir}} Y$.

\begin{example} [\textbf{Metric on birational classes}]
Suppose that $W$ is another smooth projective variety of dimension $n$. By taking fibre products over $W$, one sees that
\[  \corrdeg(X, Y) \ \le \ \corrdeg(X,W)\cdot \corrdeg(W,Y). \]
This means that $\log \big( \corrdeg \big)$ defines a metric on the  birational equivalence classes of varieties of fixed dimension.
\end{example}

\begin{example}[\textbf{Relation with degree of irrationality}] \label{Comparison.with.Deg.Irrat}
Given a projective $n$-fold $X$, recall that the \textit{degree of irrationality} $\irr(X)$ is defined to be the least degree of a rational covering $X \dra \PP^n$. This invariant has been the focus of considerable attention in recent years (\cite{BCD}, \cite{BDELU}, \cite{Chen.Stapleton}). Given generically finite coverings 
\[ X \lra \PP^n \ \ , \ \ Y \lra \PP^n\]   with $X$ and $Y$ irreducible, one can arrange via post-composing with a generic automorphism of $\PP^n$ that
 $X  \times_{\PP^n} Y $ is likewise irreducible (\cite[3.3.10]{PAG}). It follows from this that
\begin{equation} \label{Corr.Deg.vs.Irrat.Eqn}
\corrdeg(X, Y) \, \le \, \irr(X) \cdot \irr(Y). 
\end{equation}
As stated in the Introduction, our intuition is that equality in \eqref{Corr.Deg.vs.Irrat.Eqn} indicates vaguely speaking that $X$ and $Y$ are ``birationally independent." \qed
\end{example}

\begin{example} [\textbf{Yang's theorem}]  It follows from the definition or the previous remark that
\begin{equation} \label{Corrdeg.w.Pn.Eqn}
\corrdeg(X, \PP^n) \, \le \, \irr(X). 
\end{equation}
Yang \cite{Yang} proves the interesting result that equality holds when $X \subseteq \PP^{n+1}$ is a very general hypersurface of degree $d \ge 2n + 2$. We suspect that in general, the inequality \eqref{Corrdeg.w.Pn.Eqn} can be strict. 
\end{example}

As in the earlier work \cite{LP}, \cite{BCD}, \cite{BDELU} on measures of irrationality,  a natural strategy for bounding this invariant is to exploit tension between the positivity properties of the canonical bundle and vanishing statements coming from Hodge theory. We review next the various inputs to this approach. For more details, the reader can consult the papers just cited.

Fix a correspondence $Z$   as in \eqref{Correspondence.Setup}.
Via the trace on holomorphic $n$-forms from $Z$ to $X$ and $Y$,  such a correspondence  determines homomorphisms:
\begin{equation} \label{Trace.Equations}
\begin{aligned}
Z_* \  = \ \textnormal{Tr}_b \circ a^* \, &: \, H^{n,0}(X) \lra H^{n,0}(Y) \\
Z^* \   = \  \textnormal{Tr}_a\circ  b^* \, &: \, H^{n,0}(Y) \lra H^{n,0}(X).\end{aligned}
\end{equation}
(For notational purposes only -- to choose between the covariant $Z_*$ and the contravariant $Z^*$ -- we break the symmetry of the setup by viewing $Z$ as a correspondence from $X$ to $Y$.) Concretely, fix a holomorphic $n$-form on $X$ and a general point of $y \in Y$. Given $x \in b^{-1}(y)$, $db_x$ identifies $T_xX$ with $T_yY$, and then
\[    Z_*\big(\omega\big)(y) \, = \,\sum_{x \in f^{-1}(y)} \omega(x). \]
In particular, if we are able to find $\omega$ vanishing at all but one of the points of $f^{-1}(y)$ and not at the remaining one, it will follow that $Z_*(\omega) \ne 0$. Analogous considerations hold for $Z^*$. It is at this point that   birational positivity of canonical  bundles enters the picture. 

Let $X$ be an irreducible projective variety of dimension $n$, and let $L$ be a line bundle on $X$. We say that $L$ \textit{birationally separates} $p + 1$ points on $X$ -- or that it satisfies \textit{Property} $\BVA{p}$ -- if the following holds:
\begin{equation}
\parbox{5in}{There exists a non-empty Zariski-open subset $U \subseteq X$ with the property that given any $p+1$ distinct points \[ \hskip -2in x_0, \ldots, x_p \in U,\] there exists a section $s \in \HH{0}{X}{L}$ that vanishes at $x_1, \ldots, x_p$ but is non-vanishing at $x_0$.}
\end{equation}
(The definition of Property $\BVA{p}$ in \cite{BDELU} required that $L$ satisfy the analogous statement for subschemes of length $p+1$, but the present simplification works just as well for our purposes.)  For example, if $X$ maps birationally to its image in a projective space $\PP$, and if $H_X$ denotes the pullback of the hyperplane divisor fo $X$, then $\OO_X(pH)$ satisfies $\BVA{p}$. More generally, we will say that a subspace  $V \subseteq \HH{0}{X}{L}$ satisfies $\BVA{p}$ if the condition is satisfied by sections in $s \in V$.  We refer to \cite[\S 1]{BDELU} for further details.

We may summarize the discussion so far in 
\begin{lemma} \label{BVA.Implies.Non.Vanishing.Lemma}
In the situation of \eqref{Correspondence.Setup}, suppose that $K_X$ $($or a subspace $V \subseteq H^{n,0}(X) )$ satisfies $\BVA{p}$. If $
\deg(b) \le p+1$, then the homomorphism
\[
Z_* \, : \, H^{n,0}(X) \lra H^{n,0}(Y) 
\]
$($or the restriction $Z_* | V )$ is non-zero. Similarly, if $K_Y$  $($or a subspace $W \subseteq H^{n,0}(Y) )$ satisfies $\BVA{p}$ and $\deg(a) \le p + 1$, then $Z^* \ne 0$ $($or $Z^*|W \ne 0)$. \qed
\end{lemma}
\noi In \cite{LP}, \cite{BCD} and \cite{BDELU}, the Lemma was mainly applied in situations where $H^{n,0}(Y) = 0$. However as was also noted in those (and other) papers, under suitable Hodge-theoretic independence assumptions, one can infer that $Z_*$ or $Z^*$ vanish even when the spaces of forms involved are non-zero.  

Specifically, a correspondence $Z$ as in \eqref{Correspondence.Setup} gives rise to morphisms of weight $n$ Hodge structures
\[    Z_* : H^n(X) \lra H^n(Y) \ \ ,\ \ Z^* : H^n(Y)\lra H^n(X), \]
and the homomorphisms in \eqref{Trace.Equations} are the $(n,0)$ components of these. There are several situations where one knows  the vanishing of such morphisms in generic settings:

\begin{enumerate}
\item[{}] If $X$ and $Y$ are very general smooth curves of genera $g_X, g_Y \ge 1$, then \begin{equation}\label{Hodge.Indpendence.Smooth.Curves}\begin{aligned} \Hom_{\textnormal{HS}}\big( H^1(X)\, , \, H^1(Y) \big) &= 0 , \\ \Hom_{\textnormal{HS}}\big( H^1(Y)\, , \, H^1(X) \big) &= 0 .\end{aligned}\end{equation}
\vskip5pt
\item[{}]  If $X \, , \, Y \subseteq  \PP^{n+1}$ are very general hypersurfaces of degrees $d, e \ge n+2$, then
\begin{equation}\label{Hodge.Independence.General.Hypersurfaces}\begin{aligned} \Hom_{\textnormal{HS}}\big( H^n(X)_\prim\, , \, H^n(Y)  \big) &= 0 , \\ \Hom_{\textnormal{HS}}\big( H^n(Y)_\prim\, , \, H^n(X) \big) &= 0. \end{aligned}\end{equation}
\end{enumerate}
(See \cite{Ciliberto.Endos}, \cite{Voisin.Book}.) Recalling that $(n,0)$ cohomology is in any event primitive, we find:
\begin{proposition} \label{Hodge.Independ.Implies.Van.Trace.Prop}
In either of the settings just described, let $Z$ be a correspondence between $X$ and $Y$. Then the  induced homomorphisms $Z_*$ and $Z^*$ in \eqref{Trace.Equations} vanish. \end{proposition}
\noi In  \S 2  we will discuss a related hypothesis that accommodates very general hypersurfaces of large degree in an arbitrary variety.

To illustrate these considerations, we conclude this subsection by spelling out what we can say so far about very general curves and hypersurfaces.

\begin{example}[\textbf{Hodge-independent  curves}] \label{Very.General.Curve.Example} Let $X$ and $Y$ denote   curves of genera $g_X, g_Y \ge 1$ that satisfy   \eqref{Hodge.Indpendence.Smooth.Curves}.  Then
\begin{equation} \label{Curve.Gonality.Eqn}
\corrdeg(X,Y) \ = \ \gon(X) \cdot \gon(Y).
\end{equation}
In fact, \eqref{Hodge.Indpendence.Smooth.Curves} is equivalent to asking that
\[   \Pic(X \times Y) \ = \ \Pic(X) \times \Pic(Y). \]
Assuming this let $\ol{Z} \subseteq X \times Y$ be the graph of a correspondence as in \eqref{Correspondence.Setup}. Then 
\[   \OO_{X\times Y}(\ol{Z}) \ = \ B \boxtimes A \]
for some  line bundles $B$ and $A$ on $X$ and $Y$, and necessarily $h^0(B), h^0(A) \ge 2$, with $B$ and $A$ globally generated . Therefore
\begin{equation} 
\label{Curve.Corresp.Ineq} \deg(Z \rightarrow X ) =  \deg(A)     \ge \gon(Y) \ \ , \ \ \deg(Z\rightarrow Y) = \deg(B) \ge \gon(X). \end{equation}
Now assume that equality holds in \eqref{Curve.Gonality.Eqn}
 and that $Z$ computes the correspondence degree. Then both inequalities in \eqref{Curve.Corresp.Ineq} must be equalities, and therefore $h^0(B) = h^0(A) = 2$. It follows that $\ol{Z}$ is defined in $X \times Y$ by an equation of the form $s_0t_1 - s_1t_0 = 0$ where $s_0, s_1$ are sections of $B$ and $t_0, t_1$ are sections of $A$. This implies that $\ol{Z}$ is the fibre product of the coverings $X \lra \PP^1$, $Y \lra \PP^1$ defined by $B$ and $A$ respectively. 
\end{example}

\begin{example} [\textbf{Very general hypersurfaces}] \label{Weak.Bound.Corr.Def.Hypsfs} Consider next  hypersurfaces
\[  
X \, = \, X_d   \ , \ Y\, = \, Y_e \ \subseteq \ \PP^{n+1} 
\]
of degrees $d , e \ge n+2$, and let $Z$ be a correspondence as in \eqref{Correspondence.Setup}. If \eqref{Hodge.Independence.General.Hypersurfaces} holds, then 
\[
\deg(a) \, \ge \, e - n \ \ , \ \ \deg(b) \, \ge \, d - n.
\]
Indeed, should either of these inequalities fail,  Proposition \ref{Hodge.Independ.Implies.Van.Trace.Prop} would contradict Lemma \ref{BVA.Implies.Non.Vanishing.Lemma}. In the next section we will show that  stronger statements hold when $X$ and $Y$ are very general and $d, e \ge 2n + 2$.  In particular, these bounds hold when $X$ and $Y$ are very general hypersurfaces. \qed
\end{example}

\subsection*{Joint covering invariants} Invariants of covering families of curves give interesting measures of irrationality \cite{Bast}, \cite{BDELU}, \cite{BCFS}. The present subsection defines analogues for pairs of varieties. 

Let $V$ be a smooth projective variety of dimension $n$. Recall from \cite[\S1]{BDELU} that a covering family of curves on $V$ consists of a smooth family 
\[ \pi : \mathcal{C} \lra T\] of irreducible projective curves parametrized by an irreducible variety $T$, together with a dominant morphism 
\[  f : \CCC \lra V \]
with the property that for a general point $t \in T$, the map
\begin{equation} \label{Covering.Fam.Eqn}
 f_t \, : \, C_t =_{\text{def}} \pi^{-1}(t) \lra V \end{equation}
is birational to its image. We say that the family has genus $g$ or  gonality $c$ if the general fibre $C_t$ has those invariants. The   \textit{covering genus} and  \textit{covering gonality}
\[   \covgenus(V) \ \ , \ \ \covgon(V) \]
of $V$ are defined to be the minimum gonality and minimum genus of covering families of curves on $V$. Evidently we can restrict attention if we like to covering families with $\dim T = n-1$, but this is not necessary.

\begin{remark} \label{Non.Birational.Familites}
Observe that if we start with a covering family  of curves
\[  \pi^\pr : \CCC^\pr \lra T^\pr  \ \ , \ \ f^\pr : \CCC^\pr \lra Z \] that fails to satisfy the birationality condition \eqref{Covering.Fam.Eqn}, we can map it to one that does by resolving the singularities of  the image of $\CCC^\pr$ in $Z \times T^\pr$. Since gonality and genus do not increase under coverings, the invariants of this new family are no larger than those of $\CCC^\pr$.
\end{remark}

Now suppose that $u : Z \lra X \times Y$ is a correspondence between $n$-dimensional varieties  $X$ and $Y$ as in \eqref{Correspondence.Setup}. We refer to a covering family of curves on $Z$ as a joint covering family of  $X$ and $Y$. Note that the general curve $C_t$ in such a family maps birationally to its image in $X \times Y$, but it might fail to map birationally to its image in $X$ or $Y$. 
\begin{definition} The \textit{joint covering genus} and  \textit{joint covering gonality} of $X$ and $Y$ are  defined to be the minimum genus or gonality of such joint families:
\begin{align*}
\covgenus(X, Y) \ &= \ \min_Z \big\{ \covgenus(Z) \big \}\\
\covgon(X, Y) \ &= \ \min_Z \big\{ \covgon(Z) \big \},
\end{align*}
the minima being taken over all correspondences $Z$ as in \eqref{Correspondence.Setup}.
\end{definition}

\begin{remark} \label{Family.of.Curves.Remark}
Let $\pi : \CCC \lra T$ be a family of curves of genus $g$ and gonality $c$  that admits a mapping $f : \CCC \lra X \times Y$  dominating each factor and that has the property that the generic member $C_t$ maps finitely to its images in $X$ and $Y$. Then $\CCC$ contributes to the computation of the joint covering invariants in the sense that \[
\covgenus(X,Y) \, \le \, g \ \ , \ \ \covgon(X,Y) \, \le \, c. 
\]
In fact, after possibly replacing $T$ by a suitable subvariety, we can suppose first that $\dim T = n-1$. Then replacing $T$ by an  open subset we can find a correspondence $Z$ between $X$ and $Y$ and a dominant mapping $\CCC \lra Z$. Keeping in mind Remark \ref{Non.Birational.Familites}, the assertion follows.
\end{remark}

\begin{example}
It follows from the definitions and Remark \ref{Non.Birational.Familites} that
\begin{align*}
\covgenus(X, Y) \ &\ge \ \max \big \{ \covgenus(X)\, , \, \covgenus(Y) \big \}\\
\covgon(X, Y) \ &\ge \ \max \big \{ \covgon(X)\, , \, \covgon(Y) \big \}.
\end{align*}
When $Y = \PP^n$, equality holds:\[
\covgenus(X, \PP^n) \,= \, \covgenus(X) \ \ , \ \ 
\covgon(X, \PP^n) \, = \, \covgon(X). \ \ \ \qed
\]
\end{example}

As in Example \ref{Comparison.with.Deg.Irrat}, these invariants satisfy multiplicative upper bounds:
\begin{proposition}[\textbf{Upper Bounds}] \label{UB.Joint.Covering.Invariants}
Let $X$ and $Y$ be smooth projective varieties of dimension $n$. Then
\begin{align}
\label{Joint.UB1} \covgon(X, Y) \ &\le \ \covgon(X) \cdot \covgon(Y) \\
\label{Joint.UB2} \covgenus(X,Y) \ &\lesssim \ 3 \cdot \covgenus(X) \cdot \covgenus(Y).
\end{align}
\end{proposition}
\noi Denoting by $g_X, g_Y$ the covering genera of $X$ and $Y$, the precise inequality in \eqref{Joint.UB2} is that 
$\covgenus(X,Y) \le 3g_Xg_Y +  (g_X + g_Y) .$

\begin{proof} [Sketch of Proof of Proposition \ref{UB.Joint.Covering.Invariants}]
Let $\CCC = \{ C_t \}_{t \in T}$ and $\DDD = \{ D_s\}_{s \in S}$ be covering families of curves on $X$ and $Y$ respectively. After perhaps shrinking $T$ and $S$, we will  specify a family of smooth irreducible curves
\[    E_{t,s} \, \subseteq \, C_t \times D_s \]
of controlled gonality or genus that dominate  $C_t$ and $D_s$. As $t\in T$ and $s \in S$ vary, the $E_{t,s}$  sweep out a family $\mathcal{E}$ of curves mapping to $X \times Y$,  and    Remark \ref{Family.of.Curves.Remark}  then   applies to give   the stated bounds. The $E_{t,s}$ in turn are to be  divisors of general sections of a line bundle on $C_t \times D_s$ of the form $A_t \boxtimes B_s$ where $A_t$ and $B_s$ are suitable basepoint-free pencils bundles on $C_t$ and $D_s$. What remains is to indicate the choice of $A_t$ and $B_s$.

For $\eqref{Joint.UB1}$, assume that $\CCC$ and $\DDD$ compute the covering gonalities of $X$ and $Y$ respectively and take $A_t$ and $B_s$ to be the minimal degree pencils on $C_t$ and $D_s$.\footnote{We allow ourselves to pass to \'etale neighborhoods $T^\pr \lra T$ and $S^\pr \lra S$ in order to guarantee that $A_t$ and $B_s$ are rationally defined.} Then 
\[  \gon(E_{s,t}) \ \le \ \deg(A_t) \cdot \deg(B_s), \]
as required. For \eqref{Joint.UB2}, suppose that $\CCC$ and $\DDD$ compute the covering genera $g_X$ and $g_Y$ of $X$ and $Y$. Now we take $A_t$ and $B_s$ to be general pencils of degree $g_X + 1$ and $g_Y + 1$ respectively. Then the adjunction formula gives
\[
2g(E_{t,s})  \ = \ 6 g_X g_Y +  2(g_X + g_Y), 
\]
and we are done. 
\end{proof}

By way of example, we conclude this subsection with the case of $K3$ surfaces. Recall  that if $S$ is such a  surface, then 
\[
\covgenus(S) \, = \, 1 \ \ , \ \ \covgon(S) \, = \, 2. 
\]
We show that the same statements hold for pairs of $K3$'s. 
\begin{proposition} [\textbf{K3 surfaces}]
Let $S_1$ and $S_2$ be projective $K3$ surfaces. Then 
\[   \covgenus(S_1, S_2) \, = \, 1 \ \ , \ \ \covgon(S_1, S_2) \, = \, 2. \]
\end{proposition}
 \begin{proof}
 We will prove more generally that if $X_1$ and $X_2$ are smooth projective surfaces covered by non-isotrivial one-dimensional families of elliptic curves, then 
 \[
 \covgenus(X_1, X_2) \, \le \, 1, \]
 and therefore $\covgon(X_1, X_2) \le 2$. The existence of such families on $K3$ surfaces is established by Chen and Gounelas \cite[Theorem A and  Corollary on p. 2]{Chen.Gounelas}. 
 
 We need to construct a family of elliptic curves that jointly covers $X_1$ and $X_2$. To this end
 let $\pi_1 : \EEE_1 \lra T_1$ and $\pi_2 : \EEE_2 \lra T_2$ be non-isotrivial families of elliptic curves   possessing dominant maps
 \[  f_1 : \EEE_1 \lra X_1 \ \ , \ \ f_2 : \EEE_2 \lra X_2.\] By adding suitable level data, we can pass to coverings $T_1^\pr \ra T_1,    T_2^\pr \ra T_2$ of $T_1$ and $T_2$ so that the resulting families
$ \EEE^\pr_1 = T_1^\pr \times_{T_1} \EEE_1$ and $ \EEE^\pr_2 = T_2^\pr \times_{T_2} \EEE_2 $
 pull back from a universal family $\EEE \lra M$ under moduli maps
$\mu_1 :   T_1^\pr \lra M $ and $ \mu_2 : T_2^\pr \lra M$.  Observe that these maps are dominant since $\dim M = 1$ and our original families are non-isotrivial. 
Now let $T$ be an irreducible component of $T_1^\pr \times _M T_2^\pr$ that dominates both $T_1^\pr$ and $T_2^\pr$:
\[
\xymatrix@R= 15pt{
&  & T \ar[dl] \ar[dr] \\
T_1 & T_1^\pr \ar[l] \ar[dr]_{\mu_1}& & T_2^\pr \ar[r] \ar[dl]^{\mu_2}& T_2. \\
& & M
}
\]
Then
\[ 
\EEE_1 \times_{T_1} T \  \cong \  \EEE \times_M T \  \cong \ \EE_2 \times_{T_2}     T \]
over $T$, and we arrive at a common family. 
 \end{proof}

\begin{remark}[\textbf{Abelian surfaces}] As in the case of $K3$'s, if $A$ is a very general abelian surface then $\covgenus(A) = 2$ and $\covgon(A) = 2$. However here the story for joint covering invariants is different: we show in Remark \ref{cov.gon(A,B)} that the joint covering gonality of two very general abelian surfaces is $\ge 3$ (and hence the joint covering genus is likewise $\ge 3$). \qed \end{remark}

\section{The correspondence degree of very general hypersurfaces}

This section contains the proof of Theorem \ref{TheoremA}
 from the Introduction. We also establish an asymptotic bound for  the correspondence degree between two very general large degree hypersurfaces  in an arbitrary variety.

 \subsection*{Proof of Theorem \ref{TheoremA}} Much of the argument  follows the line of reasoning in \cite{BDELU}, so we will  be brief.
  
 Fix very general hypersurfaces
\[  X \, =\, X_d \ \ , \ \ Y \, = \, Y_e \ \subseteq \ \PP^{n+1} \]
of degrees $d, e \ge 2n+2$, and let $Z$ be a correspondence between $X$ and $Y$ as in \eqref{Correspondence.Setup}:
\[
\xymatrix{
& Z \ar[dl]_a\ar[dr]^b  \\
X & & Y.
 }
\]
 We may -- and do -- assume that $X$ and $Y$ satisfy \eqref{Hodge.Independence.General.Hypersurfaces}. Therefore both of the homomorphisms
\[ Z_* \, : \, H^{n,0}(X) \lra H^{n,0}(Y)  \ \ , \ \ Z^* \, : \, H^{n,0}(Y) \lra H^{n,0}(X) \]
vanish: in the terminology of \cite{BCD}, this means that $Z$ is a \textit{null-trace correspondence} in both directions. We need to show that
\[   \deg(b) \, \ge \, d -1 \ \ , \ \ \deg(a) \, \ge \, e-1,\]
and thanks to Example \ref{Weak.Bound.Corr.Def.Hypsfs} we already know that 
$\deg(b)  \ge d-n $\ and $\deg(a)   \ge   e-n$. The situation being symmetric, we focus on the mapping  $b$. Set
\[ \delta \  =_\text{def}\ \ \deg(b). \]  We assume that $\delta \le d-1$, aiming for a contradiction when $\delta \le d-2$ and an analysis of the geometry when $\delta = d-1$. 

Fix next a general point $y \in Y$. The fibre of $Z$ over $y$  sits naturally as a subset of $X$ and hence also $\PP^{n+1}$:
\[ Z_y \, =_{\text{def}} \, b^{-1}(y)\, \subseteq \, X\, \subseteq \, \PP^{n+1}.\]  Since $\delta \le 2d - 2n + 1$, it follows from \cite[Theorem 2.5]{BCD} and the vanishing of $Z_*$ that
\begin{equation}
\text{The finite set $Z_y$ spans a {line} $\, \ell_y \, \subseteq \,\PP^{n+1}$.}
\end{equation}
(In brief, $Z_y\subseteq \PP^{n+1}$ consists of $d-n\le \delta\le d-1$ points that are not separated by $\OO_{\PP}(d-n-2)$.  Bastianelli-Cortini-DePoi show that this forces $Z_y$ to lie on a line.)

As in \cite{BDELU} and \cite{Yang}, the idea of the proof is quite simple. Write
\[     X   \cdot  \ell_y \ = \ Z_y \, + \, F_y, \]
where $F_y$ is a zero-cycle of degree $d - \delta$. Letting $y$ vary in $Y$, the $F_y$ (or suitable subcycles thereof) sweep out a subvariety $S \subseteq X$ of dimension $0 \le s \le n-1$. The  new ingredient is a cohomological argument  showing that if $ s \ge 1$, then  $S$ is covered by zero-cycles of degree $\le n$ that fail to impose independent conditions on  $H^{s,0}(S)$.\footnote{In the setting of \cite{BDELU} and \cite{Yang}  this followed for free using the rationality of $Y$.} The proof then concludes as in \cite{BDELU}: computations of Ein \cite{Ein} and Voisin \cite{Voisin} show that a very general hypersurface $X$ does not contain such a subvariety.

Turning to details, denote by  $\GG = \GG(1, n+1)$ the Grassmannian parameterizing lines in $\PP^{n+1}$, and let $Y^\pr$ be a smooth variety birational to $Y$ on which the rational mapping
\[    Y \dra \GG \ \ , \ \ y \mapsto \ell_y \]
resolves to a morphism $Y^\pr \lra \GG$. The point-line incidence correspondence over $\GG$ pulls back to a $\PP^1$-bundle $\pi : W^\pr \lra Y^\pr$ that admits a surjective mapping $\mu : W^\pr \lra \PP^{n+1}$. Set
$Z^* \ = \ \mu^*(X) \subseteq W^\pr$. Then we can write
\[ Z^* \ = \ Z^\pr \, + \, F, \]
where $Z^\pr \subseteq W^\pr$ is a reduced and irreducible divisor   of relative degree $\delta$ over $Y^\pr$, whose fibre over a general point $y^\pr \in Y^\pr$ is $Z_{y^\pr}$, and $F$ is effective of relative degree $d - \delta \le n$ over $Y^\pr$. Note that $Z^\pr$ is birational to $Z$, and that  it (or, strictly speaking, a desingularization) represents birationally the original correspondence $Z$.

Fix an irreducible component $V_0$ of $F$ that dominates $Y^\pr$, and set
\[   S \, =_{\text{def}} \, \mu(V_0) \, \subseteq \, X \ \ , \ \ s \, =_\text{def} \, \dim S \ \ , \ \ e \, = \, \deg(V_0/Y^\pr). \]  
Evidently $s \le n-1$ since $S \subseteq X$. 
Assuming for the time being that
$s \ne 0$
 we will derive a contradiction.  Denote by $V_1$ the image of $V_0$ in $S \times Y^\pr$: then
\[   e^\pr   =_{\text{def}} \deg(V_1\rightarrow Y^\pr) \, \big | \, e, \] and the general fibre of $V_1$ over $Y$ sits as a subset of $S$. After replacing $S$ and $V_1$ by suitable desingularizations, we arrive at a diagram
\[
\xymatrix@C = 2pc @R=20pt{
X & S^\pr \ar[l]& V^\pr \ar[l]_{\mu_{V^\pr}}  \ar[d]^{\pi_{V^\pr}} \\ & & \, Y^\pr.}
\]
As in \cite[page 2383]{BDELU}, a computation with relative canonical bundles shows that
\begin{equation} \label{p2383.Equation}
e(n-s) \, \le \, n. 
\end{equation} 

  Continuing to assume that $1 \le s \le n-1$, let $T^\pr \subseteq Y^\pr$ be a  general complete intersection of $n-s$ very ample divisors in $Y$, so that $\dim T^\pr = s$, and let 
\[ E^\pr \ = \ \pi_{V^\pr}^{*}(T^\pr) \ \subseteq V^\pr. \]
be the inverse image of $T^\pr$ in $V^\pr$. By choosing $T^\pr$ suitably we may suppose that $E^\pr$ is non-singular, that $E^\pr \lra S^\pr$ is generically finite, and that we have a commutative diagram:
\begin{equation}\label{Fibre.Square.Computation}
\begin{gathered}\xymatrix@C = 1 pc@R=.6pc{
S^\pr& &  \,V^\pr\ \ar[dd] \ar[ll] & & \ E^\pr \ar@{_{(}->}[ll]_i \ar[dd]  \\
  & &      &  \square    \\
& &  \,  Y^\pr\   & &T^\pr  \ar@{_{(}->}[ll]^j } 
\end{gathered}.
\end{equation}
Viewing $E^\pr$ as a correspondence from $S^\pr$ to $T^\pr$, we assert
\begin{italicclaim} \label{Clam.Vanishing.Hso}
The homomorphism
\begin{equation} \label{Trace.Vanishes}  H^{s,0}(S^\pr) \lra H^{s,0}(T^\pr) \end{equation}
determined by $E^\pr$ vanishes.
\end{italicclaim}
\noi Granting \ref{Clam.Vanishing.Hso}, it follows that 
$(s,0)$-forms on $S^\pr$ do not separate the $e^\pr \le e$ points in a generic fibre of $E^\pr \lra T^\pr$. In particular, $K_{S^\pr}$ does not satisfy $\BVA{e-1}$. On the other hand,  computations of Ein \cite{Ein} and Voisin \cite{Voisin}  show that if $X\subseteq \PP^{n+1}$ is a very general hypersurface of degree $d$, and if $S^\pr$ is the desingularization of a subvariety of $X$ having dimension $s \ge 1$, then $K_{S^\pr}$ satisfies $\BVA{d+s-2n-2}$. Thus $e \ge d + s -2n$, and combined with \eqref{p2383.Equation} this contradicts the assumption that $d \ge 2n+2$ (compare \cite[pp.\ 2383--2384]{BDELU}).

As for the Claim, it follows from the projection formula that in the situation of \eqref{Fibre.Square.Computation} the two morphisms
\[   V^\pr_* : \HH{s}{S^\pr}{\CC} \lra \HH{s}{Y^\pr}{\CC} \ \  , \ \ E^\pr_* :  \HH{s}{S^\pr}{\CC} \lra \HH{s}{T^\pr}{\CC}\]
of Hodge structures are related in the natural way:
\[     E^\pr_* \ = \ j^* \circ V^\pr_*. \]
Therefore the homomorphism in \eqref{Trace.Vanishes} factors through $H^{s,0}(Y^\pr)$. But $Y^\pr$ is birational to a hypersurface $Y \subseteq \PP^{n+1}$ of dimension $n > s$, and hence $H^{s,0}(Y^\pr) = H^{s,0}(Y) = 0$. 
 
The only remaining possibility is that $s = 0$. Then  $S$ consists of a single  point  $o \in X$ containing all the lines $\ell_y$. It follows  first of all that $\delta = d-1$, completing the proof that $\deg(b) \ge d-1$  (and hence by symmetry that  $\deg(a) \ge e-1$). 
Moreover  each $o \ne x \in X$ lies on a unique line through $o$,    so writing $\Lambda = \PP^{n}$ for the $n$-dimensional projective space parametrizing lines through $o$, we get  a commutative diagram of rational maps:
\begin{equation} \label{s=0.Equation}
\begin{gathered}
\xymatrix@R =15pt@C=10pt{
& Z^\pr \ar[dl]_a \ar[dr]^b \\
X\ar@{-->}[dr]_(.35){\phi}& & Y^\pr\ar@{-->}[dl]^(.35){\psi} \\
& \Lambda = \PP^n,
} 
\end{gathered}
\end{equation}
\noi where $\deg(\phi) = d-1$ and $\psi$ maps $y\in Y$ to the parameter point of the line it determines.  Note that $\deg(\psi) \ge e-1$ thanks to \cite{BDELU}. Now suppose that $Z$ has minimal degrees   $e-1$ and $d - 1$ over  $X$ and $Y$. This implies that  \eqref{s=0.Equation} is birationally equivalent to a Cartesian square. In other words, correspondences witnessing the equality $\corrdeg(X,Y) = \irr(X) \cdot \irr(Y)$ are birationally fibre products of minimal degree coverings of $\PP^n$ by $X$ and $Y$.  This establishes the last assertion of Theorem \ref{TheoremA}.

\subsection*{Hypersurfaces in an arbitrary variety}  Let $M$ be a smooth projective variety of dimension $n + 1$,  let $A$ be a very ample divisor on $M$, and let $X = X_d \in \linser{dA}$ be a smooth divisor. It was noted in \cite[Corollary 1.12, Remark 3.2]{BDELU} that
\[    \irr(X_d) \to \infty \ \ \text{as} \ \ d \to \infty. \]
We would like to prove an analogous statement for the correspondence degree between very general hypersurfaces $X_d \in \linser{dA}$ and $Y_e \in \linser{eA}$. However one cannot expect the analogue of \eqref{Hodge.Independence.General.Hypersurfaces} to hold since $M$ may produce  non-vanishing primitive cohomology  in both $H^n(X)$ and $H^n(Y)$. Therefore we work instead with vanishing cohomology. The resulting bound appears as Theorem \ref{Corr.Deg.General.X.Y.Bound}. 

Continuing to assume that $A$ is a very ample divisor on $M$, let $D \in \linser{A}$ be a smooth divisor, so that $\dim D = n$.  Recall that the \textit{vanishing cohomology} of $D$ is defined to be the kernel of the Gysin map
\[
H^n(D)_\van \ =_\text{def} \ \ker \big( \, j_* : H^n(D) \lra H^{n+2}(M) \, \big )
\]
determined by the inclusion $j : D \hookrightarrow M$. (This is the space spanned by the vanishing cycles on $D$ as one moves it in a Lefschetz pencil.) One has an orthogonal decomposition of Hodge structures
\begin{equation} \label{Vanishing.Cohom.Eqn}
H^n(D) \ = \ H^n(D)_\van \, \oplus \, j^* H^n(M). 
\end{equation}
We refer to \cite[Chapter 2.3]{Voisin.Book} for a nice discussion of these matters.

An important point for our purposes is that adjoint divisors contribute to $H^n(D)_\van$:
\begin{lemma} \label{Adjoint.is.n.0.van}
Denote by $\Adj_M(D)$ the image of the canonical mapping:
\[
\HH{0}{M}{\OO_M(K_M + D)} \lra \HH{0}{D}{\OO_D(K_D)} \, = \, H^{n,0}(D). \tag{*}
\]
Then  $\Adj_M(D)$ is the $(n,0)$ piece of the vanishing cohomology of $D$:
\[ \Adj_M(D)  \ = \ H^{n,0}(D)_\van . \]
\end{lemma}
\begin{proof}
Poincar\'e residue determines a short exact sequence
\[
0 \lra \Omega^{n+1}_M \lra \Omega^{n+1}_M(\log D) \lra \Omega^n_D \lra 0 \] 
of sheaves on $M$, and (*) comes from the corresponding long exact sequence on cohomology. On the other hand, the connecting homomorphism \[ \HH{0}{D}{\Omega^n_D} \lra \HH{1}{M}{ \Omega^{n+1}_M}\] is a Hodge component of the Gysin mapping
$
j_* : \HH{n}{D}{\CC}\lra \HH{n+2}{M}{\CC}. 
$
The assertion follows. 
\end{proof}

\begin{lemma} \label{Hodge.Indep.Van.Cohom.Lemma}
Continue to assume that $A$ is a very ample divisor on $M$.
If  $k$ is sufficiently large, then for any fixed Hodge structure $H$ of weight $n$, a very general divisor  
\[ X   \in \linser{kA}\] has the property that
 $\Hom_{\HS} \big ( H^n(X)_\van ,H  \big)= 0$. 
 \end{lemma}
\begin{proof}[Sketch of Proof] We can assume without loss of generality that $H$ is irreducible. A standard argument using the irreducibility of the monodromy action on the vanishing cohomology of divisors in a Lefschetz pencil (\cite[\S 3.2.3]{Voisin.Book}) shows that as soon as $k \ge 1$, the Hodge structure on $H^n(X)_\van$ for very general $X$ is either irreducible or a direct sum of copies of a single irreducible Hodge structure. Hence if the conclusion of the Lemma fails, $H^n(X)_\van$  is isomorphic to a  sum of copies of $H$ for generic $X \in \linser{kA}$. By \eqref{Vanishing.Cohom.Eqn}, this implies that $H^n(X)$  is   itself a fixed Hodge structure for very general $X$. But this contradicts Green's theorem \cite{Green} that local Torelli holds for divisors in $\linser{kA}$ provided that $k$ is sufficiently large. 
\end{proof}

 \begin{theorem}  \label{Corr.Deg.General.X.Y.Bound}
Let $A$ be a very ample divisor on $M$, and suppose that $A$ satisfies $\BVA{a}$. There exists an integer $b = b(M,A) \ge 0$ $($depending on $M$ and $A)$ with the property that if 
\[
X  = X_{d+b} \in \ \linser{(d+b)A} \ \ , \ \ Y = Y_{e+b} \in \linser{(e+b)A} 
\] 
are  very general divisors, then
\begin{equation} \label{Lower.Bound.Corr.Deg.Eqn} \corrdeg(X,Y) \ > \ (da)\cdot (ea) .  \end{equation}
\end{theorem}
\begin{proof}
To begin with, choose $b =b(M,A)$ so that $bA + K_M$ is effective and so that the conclusion of Lemma \ref{Hodge.Indep.Van.Cohom.Lemma}
 holds for every $k \ge b$. We may  then suppose that  \[
 \Hom_{\HS}\big(H^n(X)_\van, H^n(Y)\big) \,  = \, 0 \ \ \text{and} \ \ \Hom_{\HS}\big(H^n(Y)_\van, H^n(X)\big) \,  = \, 0.
 \]
Now fix a correspondence $Z$ between $X$ and $Y$. In view of Lemma \ref{Adjoint.is.n.0.van}, it follows that the homomorphisms
 \[
 Z_*  :  \Adj_M(X)   \lra H^{n,0}(Y) \ \ , \ \ Z^*  :  \Adj_M(Y) \lra H^{n,0}(X)
 \]
 arising from the restrictions ot $Z_*$ and $Z^*$ to vanishing cohomology are the zero mappings. On the other hand, it is elementary that if  $\BVA{p}$  holds for a line bundle $L$ on a variety,  then $mL$ satisfies $\BVA{mp}$ for any $m \ge 1$. In our situation, this means that  $K_M + (m+b)A$ satisfies $\BVA{ma}$, and therefore $\Adj_M(X)$ birationally separates $da+1$ points on a general choice of $X$ while $\Adj_M(Y)$ birationally separates $ea+1$ points on general $Y$. Thus
 \[
 \deg(Z \rightarrow X) \, > \, ea \ \ , \ \  \deg(Z \rightarrow Y) \, > \, da
 \]
 thanks to Lemma \ref{BVA.Implies.Non.Vanishing.Lemma}, and the result follows. 
\end{proof}

\begin{example} [\textbf{Asymptotic sharpness of Theorem \ref{Corr.Deg.General.X.Y.Bound}}] \label{Product.CorrDeg.Example} A simple construction shows that \eqref{Lower.Bound.Corr.Deg.Eqn} can be asymptotically sharp. Specifically, suppose that $M = N \times \PP^1$ for some smooth $n$-fold $N$ and that $A$ has degree one on the fibres of $M \rightarrow N$, so that $a = 1$. Taking fibre products over $N$ shows that then
\[ \corrdeg(X , Y ) \, \le \, (d+b)(e+b) \, \approx \, de  \ \ \qed \]
\end{example}

\begin{remark} [\textbf{Comparison with degrees of irrationality}] \label{Example.Compare.with.Degree.of.Irrat}
In view of Example \ref{Comparison.with.Deg.Irrat}, it is interesting in the setting of the Theorem to compare $\corrdeg(X,Y)$
with the product of the degrees of irrationality of $X_{d+b} $ and $Y _{e+b}$. While one can hope for  stronger statements, in any event 
\[   \irr(X) \, \le \, \deg_A(X) \, = \, (d+b) \cdot \deg_A(M) \]
and similarly for $Y$. Therefore  
\[
\frac{ \corrdeg( X, Y) }{ \irr(X) \cdot \irr(Y) } \ > \ \frac{a^2}{4 \deg_A(M)^2}
\]
 provided that $d, e \ge b(M,A)$. In particular, this quantity is bounded away from $0$ when $d, e \gg 0$. \qed \end{remark}

 \section{Joint covering invariants, I: in the spirit of Ein and Voisin}
 
 In this section we adapt the approach of Ein   and Voisin to study joint covering invariants. 
 The idea of those authors is to combine infinitesimal calculations with considerations of positivity to analyze subvarieties of a generic member of a family. In the setting of \eqref{Correspondence.Setup}, the plan is to move $X$ and $Y$   and consider  a  correspondence $Z$ that deforms with them. For our applications $X$ and $Y$ will  be either  curves or hypersurfaces, but we  introduce the setup agnostically. 
 
 \subsection*{Joint families}
 Suppose   given smooth projective morphisms
 \[ \mu : \XXX \lra M \ \ , \ \ \nu : \YYY \lra N \]
of relative dimension $n$, where $M$ and $N$ are smooth varieties of dimensions $p$ and $q$ respectively.  We write $X_t = \mu^{-1}(t)$ and $Y_s = \nu^{-1}(s)$ for  fibres of $\mu$ and $\nu$. The first point is:
\begin{proposition} [\textbf{Ein-Voisin for joint families}]  \label{Covering.Family.Deformation.Proposition}
Suppose that $Z = Z_{t,s}$ is a correspondence as in \eqref{Correspondence.Setup}  between very general members $X = X_t$ and $Y = Y_s$ of the two families $\mu$ and $\nu$. Then for every pair of integers $i, j \ge 0$ with $i + j = n$, there is a generically surjective morphism
\[
a^* ( \Omega^{p+i}_\XXX\,  |\, X) \otimes b^* ( \Omega^{q+j}_\YYY\, |\, Y) \lra \Omega^{n}_Z \, = \, \OO_Z(K_Z)
\]
of vector bundles on $Z$. 
\end{proposition}

\begin{proof}
By a standard argument, such a correspondence $Z$ must be the fibre of a family of correspondences. More precisely,   $Z$ is realized as a very general fibre of 
a smooth projective family $\lambda: \ZZZ \lra U$ sitting in a diagram as follows:
\begin{equation} \label{Family.Diagram}
\begin{gathered}
\xymatrix@C=15pt@R=20pt{
  \ZZZ    \ar[rr]^{(\alpha_U, \beta_U)\ \ \ \  }  \ar[dr]_\lambda& & \XXX_U \times_U \YYY_U \ar[dl] \\
  &\   \  \  \ U.   \ \  } \end{gathered}
  \end{equation}
Here $U \lra M \times N$ is an \'etale morphism, and $\XXX_U$, $\YYY_U$ are the pullbacks of $\XXX$ and $\YYY$ under the projections $\rho_M: U \lra M$ and $\rho_N: U \lra N$. Moreover,  the general fibre  $Z_u = \lambda^{-1}(u)$, viewed as a correspondence between $X_{\rho_M(u)}$ and $Y_{\rho_N(u)}$, has the properties laid out in \eqref{Correspondence.Setup}.

 Write $\alpha : \ZZZ \lra \XXX$, $\beta : \ZZZ \lra \YYY$ for the natural maps. Noting that $\XXX_U \times_U \YYY_U \lra \XXX\times \YYY$ is \'etale, the derivative of the mapping $(\alpha_U, \beta_U)$ in \eqref{Family.Diagram} determines a homomorphism
\[
\alpha^* \Omega^1_\XXX \, \oplus \, \beta^* \Omega^1_\YYY \,  \lra \, \Omega^1_\ZZZ
\]
of vector bundles on $\ZZZ$. Taking $\Lambda^{p+q+n}$ of both sides yields
\begin{equation} \label{Tensor.Product.Equation}
\bigoplus_{i+j = n} \, \alpha^*  \Omega^{p+i}_\XXX \, \otimes \,  \beta^*  \Omega^{q+j}_\YYY\ \lra \ \Omega^{p+q+n}_\ZZZ = \OO_\ZZZ(K_\ZZZ). \end{equation} 
We assert:
\begin{quote}
\vskip-25pt
\begin{claim} \label{Claim.in.Proof.of.Derivative.Surjective}
  For any $i, j \ge 0$ with $i + j = n$, the mapping
\[
\alpha^*  \Omega^{p+i}_\XXX \, \otimes \,  \beta^*  \Omega^{q+j}_\YYY\ \lra \ \Omega^{p+q+n}_\ZZZ 
\]
appearing as a summand in \eqref{Tensor.Product.Equation} is a generically surjective morphism of  bundles on $\ZZZ$. 
\end{claim}
\vskip -20pt
\end{quote}
Once this is known, the  Proposition follows from the observation that $\Omega^{p+q+n}_\ZZZ |Z = \Omega^{n}_Z$. 
 As for the Claim, return to the  diagram \eqref{Correspondence.Setup} showing the restriction of \eqref{Family.Diagram}  to   general fibres, and fix a point $z \in Z$ at which both maps $a$ and $b$ are \'etale. A local calculation shows that the homomorphism of bundles
\[
a^* \Omega_X^i \otimes b^*\Omega^j_Y \lra \Omega^n_Z = \OO_Z(K_Z)
\]
is surjective at $z$, and hence generically surjective as a map of bundles on Z. The assertion of the Claim is deduced from this. \end{proof}

\subsection*{Hypersurfaces}
In this subsection we study covering  
invariants of  hypersurfaces. The first point is to apply Proposition \ref{Covering.Family.Deformation.Proposition} to prove the following:

  \begin{proposition} \label{Canonoical.Bund.Corresp.Very.Gen.Hypsfs}
Consider very general hypersurfaces
\[ X  \, , \, Y    \ \subseteq \ \PP^{n+1} \]
  of degrees $d$ and $e$, and let $Z$ be a correspondence between them as in \eqref{Correspondence.Setup}. Write $H_X, H_Y$ for the pullbacks to $Z$ of the hyperplane classes on $X$ and $Y$, and fix integers $i, j \ge 0$ such that $i + j = n$. Then
\[
K_Z \ \succcurlyeq \  \big(  d + i  - (2n + 2) \big) H_X \, + \, \big(  e + j - (2n + 2) \big)H_Y. \footnote{We write $D \succcurlyeq D^\pr$ to indicate that the difference $D - D^\pr$ of two divisors (or divisor classes)  is effective. }
\]
\end{proposition}
 \begin{proof} [Sketch of Proof]
 This follows immediately from Proposition \ref{Covering.Family.Deformation.Proposition} using the arguments and computations of Voisin \cite{Voisin}. In brief,  we take $X$ and $Y$ to be very general fibres of  the universal families
 \[
\mu : \XXX \lra M \ \ , \ \ \nu : \YYY \lra N
 \]
 of hypersurfaces of degrees $d$ and $e$. Voisin \cite[]{Voisin}  shows that the restrictions to $X$ and $Y$ of the tangent  bundles to $\XXX$ and $\YYY$ have the property that 
 \[     T_\XXX | X \otimes \OO_X(1) \ \ \text{ and } T_\YYY | Y \otimes \OO_Y(1)  \]
 are globally generated. Noting that
 \begin{align*} \Omega^{p+i}_\XXX | X \ &= \ \Lambda^{n-i} \big( T_\XXX |X \otimes \OO_X(1) \big) \otimes \OO_X\big(d + i - (2n+2)\big) \\
 \Omega^{q+j}_\YYY | Y \ &= \ \Lambda^{n-j} \big( T_\YYY |Y \otimes \OO_Y(1) \big) \otimes \OO_Y\big(e +j - (2n+2)\big),
 \end{align*}
 the assertion follows from \ref{Covering.Family.Deformation.Proposition}.
  \end{proof}

We   use this calculation to deduce:
\begin{theorem}[\textbf{Joint covering invariants for hypersurfaces}] \label{Joint.Cover.Invariants.Hypsfs}
Let $X , Y \subseteq \PP^{n+1}$ be very general hypersurfaces of degrees $d$ and $e$. Then
\begin{gather}
   \covgon(X, Y) \ \ge \  d + e - (3n + 2).  \label{joint.gon.bound.eqn}\\ 
   \covgenus(X, Y) \  \ge \ \frac{d^2 + e^2 - 3n^2}{2} \, + \, \textnormal{(lower order terms)}. \label{joint.genus.bound.eqn}
   \end{gather}
\end{theorem}

\begin{remark} [\textbf{Additivity of bounds}] \label{Bounds.Are.Additive.Remark}
The  bounds appearing in the statement are essentially additive in the invariants of the individual hypersurfaces. Indeed, for  fixed $n$ and $d \gg 0$, one has
\[ \covgon(X_d) \, \sim\,  d \ \ , \ \ \covgenus(X_d) \,  \sim \, d^2 / 2,\]
and similarly for $Y$. 
For the gonality this follows e.g.\ from \cite{BCD}, \cite{BDELU} or \cite{BCFS}, while the covering genus is computed by a simple argument that we explain in Remark \ref{Cov.Genus.Single.Hypsf}. As to upper bounds, it follows from Proposition \ref{UB.Joint.Covering.Invariants}  that
  $\covgon(X, Y) \ \lesssim \ de$, and we conjecture that the joint gonality should actually grow multiplicatively. Since $X$ and $Y$ are covered by plane curves, $\covgenus(X,Y)$ is bounded above by the genus of a divisor of type $(1,1)$ on the product of two plane curves.   This leads to 
  \[
  \covgenus(X, Y) \ \lesssim \ \frac{de(d+e)}{2},
  \]
and again one might expect the actual value to have approximately this shape. \qed
\end{remark}

\begin{proof} [Proof of Theorem \ref{Joint.Cover.Invariants.Hypsfs}]
Returning to the situation of Proposition \ref{Canonoical.Bund.Corresp.Very.Gen.Hypsfs}, consider  a family of curves
\[
\pi : \CCC \lra T \ \ , \ \ f : \CCC \lra Z 
\]
covering $Z$, and write $C$ for a general fibre. We assume that $C$ maps birationally to its image in $X \times Y$, and the issue is to bound the gonality and genus of $C$.  Denote by $h_X, h_Y$ the pullbacks to $C$ of hyperplane divisors on $X$ and $Y$, and fix $i , j \ge 0$ with $i + j = n$.    Since $K_\CCC \succcurlyeq  f^* K_\ZZZ$,  
and since $K_\CCC | C = K_C$, Proposition \ref{Canonoical.Bund.Corresp.Very.Gen.Hypsfs} shows that 
\begin{equation} \label{lowerbound.KC}
 K_C \succcurlyeq  \big(  d +i - (2n+2) \big) h_X \, + \, \big(  e +j - (2n + 2) \big)h_Y. 
\end{equation}
Statement \eqref{joint.gon.bound.eqn}  would follow immediately if we knew that $h_X$ and $h_Y$ satisfied $\BVA{1}$, but this might not be the case if $C$ doesn't  map birationally to its image in $X$ and $Y$. We will prove an essentially combinatorial lemma to circumvent this problem.

To this end, denote by
$D_1 \lra X$  and $D_2 \lra Y$
the normalizations of the images of $C$ in $X$ and $Y$.  Thus $D_1$ and $D_2$ map birationally to their images in $\PP^{n+1}$, and  $C$ maps birationally to its image in  $D_1 \times D_2$. Observe that
\[
\gon(D_1) \ge d - n \ \ , \ \ \gon(D_2) \ge e -n 
\]
since the covering gonalities of $X$ and $Y$ satisfy these inequalities. Denote by
\[   L_1 \, = \, \OO_{D_1}(d+i -  2n-2) \ \ , \ \ L_2 \, = \, \OO_{D_2}(e +j -2n -2) \]
the pullbacks to $D_1$ and $D_2$  of the indicated line bundles on $\PP^{n+1}$.  Then $L_1$ separates $d+i -  2n-1$ points on an open subset of $D_1$, and $L_2$ separates $e +j - 2n - 1$ points on an open set. Now suppose for a contradiction that $C$ admits a covering  $\phi : C \lra \PP^1$ of degree  $< d + e - 3n -2$. Lemma \ref{Bicovering.Lemma}
 below asserts that  then $L_1 \boxtimes L_2$ separates the points in a general fibre of $\phi$. But then it follows from \eqref{lowerbound.KC} that $K_C$ separates the points in a general fibre, which is not the case.  This proves \eqref{joint.gon.bound.eqn}.
 
Turning to \eqref{joint.genus.bound.eqn},  equation  \eqref{lowerbound.KC} gives a  lower bound on the genus $g(C)$ in terms of the degrees of $h_X$ and $h_Y$. Contenting ourselves with relatively rough estimates, we claim 
\[   \deg h_X \, \ge \, d-n \ \ , \ \ \deg h_Y \, \ge \, e -n. \]
Indeed, consider as above the curve $D_1$ arising as (the normalization of)  the image of of $C$ in $X$. Then 
\[
d - n \, \le \, \gon(D_1) \, \le \, \deg \OO_{D_1}(1) \, \le \, \deg h_X, 
\]
 and similarly for $h_Y$. Then  \eqref{joint.genus.bound.eqn} follows by taking $i \approx j\approx n/2$ in \eqref{lowerbound.KC}.
\end{proof}

\begin{lemma} \label{Bicovering.Lemma}
Let  $D_1, D_2$ and $C$ be smooth projective curves, with $C$ admitting a map
\[  C \lra D_1 \times D_2\]
that is birational onto its image and dominates each factor. Let $L_i$ be a line bundle on $D_i$ that separates a certain number $d_i \ge 1$ of points on a Zariski-open subset, with
\[   d_i \, \le \, \gon(D_i). \]
Finally, suppose given a covering $\phi: C \lra \PP^1$ with\[ d \, =  \, \deg(\phi)   \ < \ d_1 + d_2. \]
Then $L_1 \boxtimes L_2$ separates the points of a general fibre of $\phi$. 
\end{lemma}

\begin{proof}
We suppose to begin with that $ d_1 \ge 2$. Denote by $C_i$ the normalization of the image of $C$ under the natural map $C \lra D_i \times \PP^1$, so that one has a factorization of $\phi$:
\[
\xymatrix@C = 50pt {
& C \ar[dr]^{\phi}\ar[d] \ar[dl]_{\pro_i}  \\
D_i &C_i \ar[l]\ar[r]^{\phi_i\ \ \ } & \PP^1 .
} 
\]
Put $e_i = \deg(\phi_i)$. Observe that $e_i | d$, and 
\[d_i \, \le \, \gon(D_i) \, \le \, \gon(C_i) \, \le \, e_i. \]
Fixing a general point $a \in \PP^1$, set $T = \phi^{-1}(a) \subseteq D_1 \times D_2$, and write
\[ T_i \, = \, \pro_i(T) \, \subseteq \, D_i. \]
Thus $\# T_i = e_i$, and each fibre of the projection $T \lra T_i$ consists of precisely $d/e_i$ points.

Now choose any point $(t_1, t_2) \in T$. It suffices to construct sections
\[   s_1   \in   \GGG{D_1}{L_1} \ \ , \ \  s_2  \in   \GGG{D_2}{L_2} \]
with the property that $s_1 \boxtimes s_2$ vanishes as every point of $T$ other than $(t_1, t_2)$. To this end, consider first the set
\[    A_1 \ = \ \big \{\, t \in T_1 \mid \, (t, t_2) \in T \, \big \}.\]
Recalling that we are assuming $d_1 \ge 2$, one has
\[
\# \, A_1 \ = \ \left( \frac{d}{e_2} \right) \  \le \ \frac{d_1 + d_2 -1 }{e_2} \ \le \ \frac{d_1d_2}{e_2} \ \le \ d_1.
\]
Expanding (if necessary) $A_1$ to a $d_1$-element subset $A_1^\sharp$
with $T \supseteq A_1^\sharp \supseteq A_1$, we take $ s_1$ to be a section of $L_1$ not vanishing at $t_1$ but vanishing at every other point of $A_1^\sharp$.

Consider next the sets:
\begin{align*}
B_2 \ &= \ \big \{\, t^\pr \in T_2 \mid \, (t, t^\pr) \in T \  \text{for } t \in T_1 -  A_1^\sharp \, \big \} \\
B_2^\sharp \ &= \ B_2 \, \cup \, \big\{ \, t^\pr \mid \ (t_1, t^\pr) \in T \, \big \}. 
\end{align*}
Then
\begin{align*}
\# \, B_2^\sharp \ &= \ \left ( \frac{d}{e_1} \right) \big( (e_1 -  d_1) +1  \big) \\
&= \ d \ - \ \frac{d(d_1 - 1) }{e_1} \\
&\le \ (d_1 + d_2 -1) \, - \, (d_1 - 1) \ = \ d_2. 
\end{align*}
We then take $s_2 \in \Gamma(D_2, L_2) $ to be a section that does not vanish at $t_2$ but vanishes at every other point of $B_2^\sharp$. 

It remains to treat the possibility that $d_1 = 1$, so that $d \le d_2$. This forces $e_2 = d_2$, and by taking $s_1$ to be non-vanishing on $T_1$, we can separate points of $T$ using sections of $L_2$. 
\end{proof}

 \begin{remark} [\textbf{Covering genus of a single hypersurface}] \label{Cov.Genus.Single.Hypsf}
 Let $X_d \subseteq \PP^{n+1}$ be a smooth hypersurface of degree $d$ and dimension $n$. While there has been considerable interest in bounding the least genus of a curve on $X$, as far as we can tell the covering genus has not received much attention. In any event, the asymptotic picture is very elementary: for fixed $n$, as a function of $d$, one has:
 \[  \covgenus(X_d)  \ \sim \   {d^2} / {2} .\] 
Evidently $\covgenus(X) \le \binom{d-1}{2}$ since $X$ is covered by plane curves of degree $d$. For the reverse inequality, one argues as in the proof of \eqref{joint.genus.bound.eqn}. Indeed, recall e.g.\ from \cite{BDELU} that if $C \rightarrow X$ is the general member of a covering family, then 
\[   K_C\succcurlyeq (d-n-2)h \]
where $h$ is the pullback of a hyperplane divisor to $C$. This implies in the first place that $\gon(C) \ge (d - n)$, and therefore  that 
\[ \deg h  \,  =  \, \textnormal{(degree of image of $C$ in $\PP^{n+1}$)}  \,  \ge \,  (d-n). \]
We then find that $2g(C) - 2 \ge (d-n)(d-n-2)$, and the assertion follows. 
 The   precise covering gonality of a very general hypersurface $X$ has been (essentially) computed by Bastianelli et al in \cite{BCFS}. It would be interesting to find a more accurate estimate of its  covering genus. \qed  \end{remark}

\subsection*{Curves}  We now turn to the case of curves, and in particular prove Theorem \ref{Curve.Genus.Bound} from the Introduction.
 
In order to apply Proposition \ref{Covering.Family.Deformation.Proposition}, one needs to control the positivity of the  restricted tangent bundle of the universal deformation of a curve. The following statement is an analogue of Proposition 1.1 in \cite{Voisin}. 
\begin{lemma} \label{Voisin.Positivity.Lemma.for.Curves}
Let $X$ be a smooth projective curve of genus $g \ge 2$, and let 
\[
\mu : \XXX \lra M\] be a local universal  deformation of $X$, so that $\dim M = 3g -3$. Then for any point $x \in X$, 
\begin{equation} \label{UX.nef}
T_\XXX|X \otimes \OO_X(x) \ \text{ is  an ample vector bundle on $X$}. 
 \end{equation}
\end{lemma}
\begin{proof}
For simplicity, write $U_X =  T_\XXX | X$. This bundle is uniquely characterized by the fact that it sits in an extension
\[
0 \lra T_X \lra U_X \lra \HH{1}{X}{T_X} \otimes_\CC \OO_X \lra 0
\]
(namely the tangent/normal bundle sequence for the inclusion $X \subseteq \XXX$) for which the connecting homomorphism 
\[ \HH{1}{X}{T_X}  = \HH{1}{X}{T_X}  \otimes \HH{0}{X}{\OO_X} \lra \HH{1}{X}{T_X} \] is the identity map. This description shows that $U_X$ is the dual of the restriction to $X$ of one of the  Picard bundles on $\Jac(X)$ studied   in \cite{EL}: $U_X \cong E_{2K}^*$  in the notation of \cite{EL},   the identification following  from Lemma 2.3 of that paper.   It was established in \cite[Proposition 2.2]{EL} that $U_X$ is stable, and we see that $U_X$ has slope \[ 
\mu( U_X) \ = \ \frac{2-2g}{3g-2} \ > \ -1. 
\]  Therefore $U_X \otimes \OO_X(x)$, being stable of positive slope, is ample.
\end{proof}

Using the Lemma, we deduce:
\begin{proposition} \label{Two.Curve.Correspondence.Prop}
Let $X$ and $Y$ be very general curves of genera $g_X, g_Y \ge 2$ and let $Z$ be a correspondence between them. As in \eqref{Correspondence.Setup}  denote by
\[
a : Z \lra X \ \ , \ \ b :Z \lra Y
\] the two projections. Then
\[   \omega_Z \ = \ a^* \omega_X \otimes b^*\omega_Y \otimes P, \]
where $P$ is a line bundle on $Z$ having the property that
\[ \deg P \, > \,\max \{ - \deg a \, , \, - \deg b \}. \]
 \end{proposition}
\noi Note that we do not assert that $ P$ or any specific twist has non-trivial global sections. 

\begin{proof} [Proof of Proposition \ref{Two.Curve.Correspondence.Prop}]Let  $\mu : \XXX \lra M$ and $ \nu : \YYY \lra N$  
be local universal families of curves of genera $g_X, g_Y \ge 2$,  so that $p = \dim M = 3g_X -3$ and $q = \dim N = 3g_Y -3$. We may consider $X$ and $Y$ to be very general fibres of $\mu$ and $\nu$.

Write $U_X = T\XXX | X$. Noting that 
 $\Omega^{p}_\XXX | X = U_X \otimes \omega_X$, Proposition \ref{Covering.Family.Deformation.Proposition} (with $i = 0, j = 1$) asserts  the existence of a generically surjective morphism
 \[
 a^* (U_X \otimes \omega_X) \otimes b^* \omega_Y \lra \omega_Z  \] 
 of vector bundles on $Z$. In other words, there is a generically surjective morphism $a^* U_X \lra P$.   It follows from the previous Lemma that if $x \in X$ is any point, then the image of \[ a^* \big(U_X \otimes  \OO_X(x)\big) \lra P \otimes a^*\OO_X(x)\] is ample. Therefore   $P \otimes a^*\OO_X(x)$ has positive degree, so $\deg P + \deg(a) > 0$.  The positivity of $P \otimes b^* \OO_Y(y)$ is similar. \end{proof}

\begin{remark}
One can improve slightly the inequality on $\deg P$ by noting that the proof of Lemma \ref{Voisin.Positivity.Lemma.for.Curves} actually shows that the $\QQ$-twisted bundles \[  U_X \negmedspace<\tfrac{2}{3} x > \ \ \text{and} \ \ U_X\negmedspace< \tfrac{2g-2}{3g-2}x>\]
on $X$ are ample and nef respectively.  However this leads to only marginal numerical  improvements in the application. \end{remark}

\begin{proof}[Proof of Theorem \ref{Curve.Genus.Bound}] Let $X$ and $Y$ be very general curves of genera $g_X, g_Y \ge 2$, and let $Z$ be a smooth curve covering both. Then (as in Example \ref{Very.General.Curve.Example}) one has
\[
a \, = \, \deg ( Z \lra X) \, \ge \, \gon(Y) \ \ , \ \ b \, = \, \deg ( Z \lra Y) \, \ge \, \gon(X).
\]
The previous Proposition then gives:
\[
2g(Z) - 2  \ \ge \ \gon(Y) \cdot  ( 2g_X - 2) \, + \, \gon(X) \cdot (2g_Y -2) \, - \, \min \big( \gon(X), \gon(Y) \big).
\]
Recalling that  
\[  \gon(X) \, = \, \left [ \frac{g_X + 3}{2} \right] \ \ , \ \ \gon(Y) \, = \, \left [ \frac{g_Y + 3}{2} \right] ,
\]
this simplifies to an inequality of the shape 
\[  g(Z) \ \ge \ g_X   g_Y \, + \, \text{(affine linear expression in $g_X , g_Y$)},\]
as asserted. The upper bound on $\covgenus(X,Y)$ comes from computing the genus of a fibre product of gonality maps $X \lra \PP^1$ and $Y \lra \PP^1$. 
\end{proof}

 \section{Joint covering invariants, II: in the spirit of Pirola}\label{Pirola}

 The main result of this section is the computation of the joint covering gonality of very general pairs of hyperelliptic curves for almost all genera. The argument is inspired by the degenerational approach of Pirola; see also \cite{Voisin.Abel.Var} and \cite{Martin}.
 
 \begin{definition}
 A subvariety $Z$ of an abelian variety $A$ \textit{generates} the abelian subvariety $B\subseteq A$ if $Z$ is contained in a translate of $B$ but is not contained in a translate of a proper abelian subvariety of $B$.
 \end{definition}
\begin{example}[\textbf{Joint covering gonality of elliptic curves}]
The joint covering gonality of two elliptic curves $X$ and $Y$ is $2$. Indeed, the set of abelian surfaces which are isogenous to $X\times Y$ is dense in the moduli of principally polarized abelian surfaces $\mathcal{A}_2$. Since the Torelli locus is an open dense subset of $\mathcal{A}_2$, there is a genus $2$ curve $C$ such that $J(C)$ is isogenous to $X\times Y$. Let $C'\subseteq X\times Y$ be the image of $C$ under such an isogeny. Since $C$ generates $J(C)$ so must $C'$ generate $X\times Y$. Thus, $C'$ is a hyperelliptic curve which dominates both $X$ and $Y$.
\end{example}

\begin{theorem}\label{covgonhyp}
Let $X$ and $Y$ be very general hyperelliptic curves of genera $  g_Y \ge   g_X  \ge 2$. Then:
\[
\covgon(X,Y) \ = \ 
\begin{cases}
3 \text{ or } 4 &\text{if } \, g_X \, = \,  g_Y = 2
\\
\ \ 4 &\text{if } \, g_Y\, \ge \,  3 \ \ \qquad.
\end{cases}
\]
\end{theorem}
\noindent Theorem \ref{covgonhyp}
 will be deduced  from the following result:
\begin{theorem}\label{hypjac}
Fix an abelian variety $A$, and let $Y$ a very general hyperelliptic curve of genus $g\geq 2$. Then any hyperelliptic curve in $A\times J(Y)$ must be contracted by the projection to one of the factors.
\end{theorem}

\begin{remark}[\textbf{Joint covering gonality of abelian surfaces}]\label{cov.gon(A,B)}
If $A,B$ are very general abelian surfaces, then
\[3 \ \leq \  \covgon(A,B) \ \leq  \ 4.\]
The upper bound is given by Proposition \ref{UB.Joint.Covering.Invariants} whereas the lower bound follows from Theorem \ref{hypjac}. Indeed a very general abelian surface dominates the Jacobian of a very general genus $2$ curve. It would be interesting to know which of the two possibilities actually occurs. \end{remark}

\noindent We will begin by proving Theorem \ref{covgonhyp} contingent on Theorem \ref{hypjac}.

\begin{proof}[Proof of Theorem \ref{covgonhyp}]
 Theorem \ref{hypjac} together with the obvious upper bound from Example \ref{Comparison.with.Deg.Irrat} gives the inequalities
 \[ 3 \ \leq \  \covgon(X,Y) \ \leq \  4.\]

The lower bound can be improved to $\textup{cov.gon}(X,Y)\geq 4$ when $g_Y\geq 3$ by the Castelnuovo-Severi bound: If $D\subseteq X\times Y$ is trigonal and dominates $Y$ then the basepoint-free $g^1_3$ on $D$ maps to a base-point free $g^1_3$ on $Y$ or a curve in a base-point free $g^2_3$ on $Y$. The Castelnuovo-Severi bound ensures that a curves of genus at least $3$ cannot have both a basepoint-free $g^1_2$ and a basepoint-free $g^1_3$. Finally, a curve of genus at least $2$ cannot admit a basepoint-free $g^2_3$.
\end{proof}

\begin{example}[\textbf{Joint covering gonality of general curves of genus 1 and 2}] The joint covering gonality of a very general genus $1$ curve $X$ and a very general genus $2$ curve $Y$ is $3$. Proposition \ref{nondeghyp} gives the lower bound $$\textup{cov.gon}(X,Y) \ \geq \ 3.$$
Hence it suffices to show that there is a genus $4$ curve $C$ which dominates both $X$ and $Y$. Let $\mathcal{R}_{Y,2}$ denote the moduli space of isomorphism classes of ramified double-covers $\pi: D\to Y$, where $D$ is smooth of genus $4$. This space has dimension $2$ and the Prym map
\begin{align*}\mathcal{P}: \ \  \ \  \ \   \mathcal{R}_{Y,2}\ \ \ \   \ & \longrightarrow   \ \ \  \mathcal{A}_2\\
\ (\pi: D\to Y)  \ &\longmapsto \ \  \mathcal{P}(\pi)\end{align*}
is quasi-finite. Observe that the surface $\mathcal{P}( \mathcal{R}_{Y,2})\subseteq \mathcal{A}_2$ meets some loci of split abelian surfaces in a curve by Proposition 5 of \cite{CP}. A point of this curve corresponds to an abelian surface isogenous to $X\times E$ for some elliptic curves $E$. We have thus obtained a smooth genus $4$ curve $D$ dominating $Y$ and such that $J(D)$ has an elliptic factor isogenous to $X$. Projecting to this factor and composing with the isogeny give a dominant map $D\to X$. Since a curve of genus $4$ has gonality at most $3$, we get the bound
$$\textup{cov.gon}(X,Y) \ \leq \  3. \ \ \ \qed$$
\end{example}

The proof of Theorem \ref{hypjac} will be obtained by reducing to the case of the product of an arbitrary abelian variety with a general abelian surface.

\begin{proposition}\label{nondeghyp}
Let $A$ be any fixed abelian variety and $B$ a very general abelian surface. Then any hyperelliptic curve in $A\times B$ must be contracted by the projection to one of the factors. In particular, any hyperelliptic curve in $A\times B$ is geometrically degenerate in the sense of Ran.
\end{proposition}

Proposition \ref{nondeghyp} will be proven using the method of \cite{Pirola} by reaching a contradiction with the following rigidity statement:
\begin{lemma}[\text{Pirola} \cite{Pirola}]\label{rigidhyp}
Hyperelliptic curves on abelian varieties are rigid up to translation.
\end{lemma}
\noi(In brief: a hyperelliptic curve in an abelian surface with a Weierstrass point at the identity gives rise to a rational curve in the corresponding Kummer surface. Rational curves in Kummer surfaces are rigid since $K3$ surfaces are not uniruled.)

We now turn to the proof of Proposition \ref{nondeghyp}. As a matter of notation, given a family of varietes $\mathcal{G} \lra S$, and a map $T \lra S$, we denote by $\mathcal{G}_T \lra T$ the pull-back family over $T$. Suppose then that for a very general abelian surface $B$, the variety $A\times B$ contains a hyperelliptic curve which is not contracted by either of the projection maps to the factors. Then we can find a local universal family of abelian surfaces $\mathcal{B}\lra S$ and a flat family of irreducible curves 
 \[\mathcal{Z} \ \subseteq \  A\times \mathcal{B} \ =_{\text{def}} \ A_S\times_S\mathcal{B}\]
 such that $\mathcal{Z}_s$ is hyperelliptic for all $s\in S$. Moreover, up to replacing $S$ with an open subset and $A$ by an abelian subvariety, we can assume that $\mathcal{Z}_s$ generates $A\times \mathcal{B}_s$ for all $s\in S$.\\

There are countably-many divisors $S_\lambda\subseteq S$ along which $\mathcal{B}_s$ is isogeneous to a product $\mathcal{B}_s\sim \mathcal{E}_s^\lambda\times \mathcal{F}_s^\lambda$, where $\lambda$ encodes the isogeny and $\mathcal{E}^\lambda\lra S_\lambda$ and $\mathcal{F}^\lambda\lra S_\lambda$ are locally complete families of elliptic curves. We let $\Lambda$ be the set that indexes such loci and for all $\lambda\in \Lambda$, we write
$$p_\lambda \ : \ A\times\mathcal{B}_{S_\lambda} \ \longrightarrow \  A\times \mathcal{E}^\lambda \times_{S_\lambda}  \mathcal{F}^\lambda \  \longrightarrow \  A\times \mathcal{E}^\lambda$$
for the composition of the isogeny with the projection map.\footnote{Note that we abuse notation to denote by $S_\lambda$ a generically finite cover of $S_\lambda$ as the projection $p_\lambda$ may only be defined after a generically finite base change.} \\

We will find some $\lambda\in \Lambda$ and a curve $C\subseteq S_\lambda$ such that:\begin{itemize}
\item $\mathcal{E}^\lambda_C$ is isotrivial with fiber $E$, i.e.\ $C$ parametrizes abelian varieties isogenous to $A\times E \times F$ for a fixed elliptic curve $E$ and some (varying) elliptic curve $F$,

\vskip 5pt
\item The projection $p_\lambda(\mathcal{Z}_s)\subseteq A\times E$ varies with $s\in C$.
\end{itemize}
This will provide the desired contradiction with Lemma \ref{rigidhyp}.\\

We will use the following lemma to ensure that the curve $p_\lambda(\mathcal{Z}_s)\subseteq A\times E$ varies with $s\in C$.
\begin{lemma}
Let $\mathcal{B}\lra C$ be a non-istrovial family of abelian varieties over a one-dimensional base $C$ and let $\mathcal{Z}\subseteq \mathcal{B}$ be a flat family of irreducible subvarieties. Suppose that $\mathcal{Z}_s$ generates $\mathcal{B}_s$ for all $s\in C$. Consider a variety $X$ and a morphism \[ f:\mathcal{B}\lra X_C = X \times C\]  such that for all $s\in C$ the map $f|_{\mathcal{Z}_s}$ is a composition
\[\mathcal{Z}_s \ \xrightarrow{ \ g_s \ }  \ \mathcal{Z}_s' \ \xrightarrow{ \ h_s \ }  \  X,\]
where
$g_s$ is the restriction of an isogeny $\eta: \mathcal{B}_s\lra A$ to $\mathcal{Z}_s$, and
  $h_s$ is birational on its image.
Then the composition
$$\mathcal{Z} \ \xrightarrow{ \ f \ } \  X_C \ \longrightarrow \   X$$
is generically finite over its image.
\end{lemma}
\begin{proof}
If this composition is not generically finite on its image there is a fixed subvariety $V\subseteq X$ such that $f(\mathcal{Z}_s)=V\subseteq X$ for all $s\in C$. We can choose desingularizations $\widetilde{\mathcal{Z}}_s$ and $\widetilde{V}$ of $\mathcal{Z}_s$ and $V$ respectively so that there is a morphism
\[\widetilde{f}_s \ :  \ \widetilde{\mathcal{Z}}_s \ \longrightarrow \  \widetilde{V}\]
identifying birationally to $f|_{\mathcal{Z}_s}: \mathcal{Z}_s\lra V$. The composition
$$\text{Pic}^0(\mathcal{B}_s) \ \xrightarrow{ \  {g_s}_* \ } \  \text{Pic}^0(A)\ \longrightarrow \   \text{Pic}^0(\widetilde{\mathcal{Z}}_s) \ \xrightarrow{ \ {\widetilde{f_s}}_* \ } \ \text{Pic}^0(\widetilde{V})$$
has finite kernel since $\widetilde{f}_s$ is birational and $g_s(\mathcal{Z}_s)$ generates $A$. Moreover, note that $\text{Pic}^0(\widetilde{V})$ does not depend on the choice of the desingularization $\widetilde{V}$. Accordingly, $\text{Pic}^0(\widetilde{V})$ is a fixed abelian variety and for each $s\in C$ the abelian variety $\mathcal{B}_s$ maps to $\text{Pic}^0(\widetilde{V})$ with finite kernel. This provides a contradiction since the locus of abelian varieties admitting a morphism with finite kernel to a fixed abelian variety is discrete in any locally complete family of abelian varieties.
\end{proof}

It therefore suffices to show the following:

\begin{lemma}\label{bir}
There is a $\lambda\in \Lambda$ such that for generic $s\in S_\lambda$ the map 
 $p_\lambda|_{\mathcal{Z}_s}: \mathcal{Z}_s\lra A\times \mathcal{E}^\lambda_s$ is the composition of the restriction of an isogeny with a map which is birational on its image.\end{lemma}

Indeed, if $\lambda$ is as in the statement of the lemma, given a generic elliptic curves $E$, we can take \[C\ =_{\text{def}} \ \{s\in S_\lambda: \mathcal{E}^\lambda_s=E\} \ \subseteq \  S_\lambda.\]
Then $\mathcal{E}^\lambda_C$ is isotrivial by construction and since $E$ is generic the map $p_\lambda|_{\mathcal{Z}_s}: \mathcal{Z}_s\lra A\times \mathcal{E}^\lambda_s$ is the composition of the restriction of an isogeny with a morphism which is birational on its image for generic $s\in C$.
\begin{proof}[Proof of Lemma \ref{bir}]
For all $s\in S$ consider
\[\mathcal{D}_s \ =_{\text{def}} \ \{z-z': z,z'\in \mathcal{Z}_s\} \ \subseteq \  A\times \mathcal{B}_s.\]
Shrinking $S$ if needed, we can assume that for all $s\in S$ the surface $\mathcal{D}_s$ does not contain $\{0\}\times E$ for any elliptic curve $E\subseteq \mathcal{B}_s$. Applying the following lemma with $A\times \mathcal{E}^\lambda_s$ instead of $A$ and $F=\mathcal{F}^\lambda_s$ finishes the argument.
\end{proof}

\begin{lemma}
Let $A$ be abelian variety, $F$ an elliptic curve, and $Z\subseteq A\times F$ an irreducible curve which generates $A\times F$. Moreover, assume that $\{0\}\times F$ is not contained in the surface
\[D \ =_{\textup{def}} \ \{z-z': z,z'\in Z\} \ \subseteq \  A\times F.\]
There is an isogeny $\eta: F\lra F'$ such that the restriction of the projection map $\pi_A: A\times F'\lra A$ to the curve $Z'=_{\textup{def}}(\textup{id}_A\times \eta)(Z)$ is birational on its image.
\end{lemma}

\begin{proof}
First, consider the action $\tau: F\lra \text{Aut}(A\times F)$ of $F$ on $A\times F$ by translation on the second factor. The subgroup
$$F_Z \ =_{\text{def}} \ \{x\in F: \tau_x(Z)=Z\} \ \subseteq \  F$$
is finite. 

Now let $\eta: F\lra F'$ be the isogeny with quotient $F_Z$ and let $Z'=_{\text{def}}(\textup{id}_A\times \eta)(Z)$ so that
$$F'_{Z'} \ =_{\text{def}} \ \{x'\in F': \tau_{x'}(Z')=Z'\}\subseteq F'$$
is trivial. Moreover, note that
$$D' \ =_{\text{def}} \ \{z-z': z,z'\in Z'\} \ \subseteq \  A\times F'$$
does not contain $\{0\}\times F'$.
We claim that the restriction of the projection $\pi_A: A\times F'\lra A$ to $Z'$ is birational on its image. Consider a generic $a\in \pi_A(Z')$ and preimages $(a,x_1), (a,x_2)\in Z'$. Then $(0,x_1-x_2)\in D'\cap \{0\}\times F$ which is discrete. Accordingly, for all $(a,x)\in Z'$ we must have $(a,x+(x_1-x_2))\in Z'$, so that $x_1-x_2\in F'_{Z'}=\{0\}$. It follows that $(a,x_1)=(a,x_2)$.
\end{proof}

This completes the proof of Proposition \ref{nondeghyp}.

\begin{remark}
The proof given above is different and sowewhat simpler than the arguments of \cite{Pirola,AP,Voisin.Abel.Var,Martin}. In this previous body of work, an involved density argument was used instead of Lemma \ref{bir}. However, we will need such a density argument in the course of the reduction of Theorem \ref{hypjac} to Proposition \ref{nondeghyp}.\end{remark}

Finally, we present the proof of Theorem \ref{hypjac} which proceeds by reduction to Proposition \ref{nondeghyp}. Let $\mathcal{Y}\lra \mathcal{H}_g$ denote a locally complete family of genus $g$ hyperelliptic curves and suppose that there is a flat family of hyperelliptic curves $\mathcal{Z}\subseteq A\times J(\mathcal{Y})$ such that for all $s\in \mathcal{H}_g$ the curve $\mathcal{Z}_s$ is not contracted by either of the projection maps to the factors. Let $\mathcal{G}\lra \mathcal{H}_g$ be the universal Grassmanian of $2$-planes in the relative tangent bundle of $J(\mathcal{Y})$ at the identity. By Proposition 4 of \cite{CP}, the following locus is dense
$$\{ \; T_{B,0}\subseteq T_{J(\mathcal{Y}_s),0}: s\in \mathcal{H}_g, B\subseteq J(\mathcal{Y}_s) \text{ an abelian surface }\} \ \subseteq \  \mathcal{G}.$$
Indeed, though Proposition 4 therein is formulated as a density statement in $\mathcal{H}_g$, the proof proceeds by obtaining a density statement at the level of the Grassmanian. Since $\mathcal{H}_{g}$ has dimension $2g-1$ and all Hecke translates of the image of $\mathcal{A}_2\times \mathcal{A}_{g-2}$ in $\mathcal{A}_{g}$ have codimension $\binom{g_2+1}{2}-\left(3+\binom{g-1}{2}\right)=2g-4$, all irreducible components of this locus have dimension at least $3$. We let $\Lambda$ be a set indexing these irreducible components and denote by $T_\lambda$ the component indexed by $\lambda\in \Lambda$.\\

It suffices to show that for sufficiently many of these components the abelian surface $B$ in $J(\mathcal{Y}_t)$ varies in a $3$-dimensional family with $t\in T_\lambda$. Indeed, we can then consider the composition of the isogeny and the projection
\[A\times J(\mathcal{Y}_t) \ \longrightarrow \  A\times D\times B \ \longrightarrow \ A\times B\]
where $D$ is an abelian $(g-2)$-fold. The image of $\mathcal{Z}_t$ is then a hyperelliptic curve $Z'\subseteq A\times B$. Since $B$ is a very general abelian surface $Z'$ must be contracted by one of the projection maps to the factors. This provides the desired contradiction.

Given $\lambda\in \Lambda$, there is a morphism from $T_\lambda$ to an appropriate moduli of abelian surfaces and we denote by $d_\lambda$ the dimension of its image. In particular, there is a $d_\lambda$-parameter family of abelian surface $\mathcal{B}^\lambda\lra T_\lambda$ and a dominant morphism\footnote{Again, strictly speaking, after a finite base change. }
$$p_\lambda: \  A\times J(\mathcal{Y}_{T_\lambda}) \ \lra \  A\times \mathcal{B}^\lambda.$$

Theorem \ref{hypjac} then follows from the following proposition, which will also require the use of Pirola's method. 
\begin{proposition}\label{d_lambda=3}The set
\[\bigcup_{\lambda\in \Lambda: d_\lambda=3} T_\lambda \ \subseteq \  \mathcal{G}\]
is dense.
\end{proposition}
\begin{proof}
The work of Colombo and Pirola in \cite{CP} shows that
\[\bigcup_{\lambda\in \Lambda}T_\lambda \ \subseteq \  \mathcal{G}\]
is dense. It therefore suffices to show that
\[\bigcup_{\lambda\in \Lambda: d_\lambda<3} T_\lambda \ \subseteq \  \mathcal{G}\]
cannot be dense. Let $\Lambda'=_{\text{def}}\{\lambda\in \Lambda: d_\lambda< 3\}$ and assume that $\bigcup_{\lambda\in \Lambda'} T_\lambda\subseteq \mathcal{G}$ is dense. We will use a specialization and projection argument analogous to the one used in the proof of Proposition \ref{nondeghyp} to reach a contradiction with Lemma \ref{rigidhyp}.\\

After a finite base change the family $\mathcal{Y}_{\mathcal{G}}\lra \mathcal{G}$ has a section which gives rise to an embedding $\mathcal{Y}_{\mathcal{G}}\lra J(\mathcal{Y}_{\mathcal{G}})$. Let $\mathcal{W}\subseteq J(\mathcal{Y}_{\mathcal{G}})$ be the resulting family of hyperelliptic curves. The fact that $p_\lambda|_{\mathcal{W}_t}$ is generically finite on its image for most $\lambda\in \Lambda'$ and $t\in T_\lambda$, allows us to find 

\begin{itemize}
\item a generically finite cover of $\mathcal{G}$, which we call $\mathcal{G}$ by abuse of notation,
\vskip5pt
\item a subset $\Lambda''\subseteq \Lambda'$ such that the following subset is dense
\[\bigcup_{\lambda\in \Lambda''} T_\lambda \ \subseteq \  \mathcal{G},\]
\item a family of curves $\mathcal{W}'\lra \mathcal{G}$,
\vskip 5pt
\item a map $$p: \ \mathcal{W} \ \lra \  \mathcal{W}'$$
which identifies birationally with the restriction of $p_\lambda$ to $\mathcal{W}_{T_\lambda}$ along $T_\lambda$ for all $\lambda\in \Lambda''\subseteq \Lambda'$.
\end{itemize}
(For more details regarding these steps the reader can consult the proof of Proposition 1.4 in \cite{Voisin.Abel.Var} or Section 3.4 of \cite{Martin}. A key element in the argument is the fact that the Gauss map of $p_\lambda(\mathcal{W}_t)\subseteq \mathcal{B}_t^\lambda$ is generically finite on its image for generic $t\in T_\lambda$ and all $\lambda\in \Lambda'$. This follows from the fact $\mathcal{W}_t$ generates $J(\mathcal{Y}_t)$ for all $t\in \mathcal{G}$.)

Now replacing $\mathcal{G}$ with an open subset if needed we can consider desingularization $\widetilde{\mathcal{W}}\lra \mathcal{G}$ and $\widetilde{\mathcal{W}}'\lra \mathcal{G}$ with smooth fibers and a map
$\widetilde{p}: \widetilde{\mathcal{W}}\lra \widetilde{\mathcal{W}}'$ which identifies birationally with $p$ on fibers. Consider the morphism
$$\textup{Pic}^0(J(\mathcal{Y}_t)) \ \lra \  \textup{Pic}^0(\widetilde{\mathcal{W}}_t) \ \xrightarrow{ \ \widetilde{p}_* \ }  \ \textup{Pic}^0(\widetilde{\mathcal{W}}_t').$$
One easily checks that for generic $t$ this morphism is non-zero. Since $\textup{Pic}^0(J(\mathcal{Y}_t))$ is simple for generic $t\in \mathcal{G}$, it follows that the composition above has finite kernel for generic $t\in \mathcal{G}$. Shrink $\mathcal{G}$ further so that the composition has finite kernel for all $t\in \mathcal{G}$.\\

Now consider $\lambda\in \Lambda''$ such that the image of $T_\lambda$ in $\mathcal{G}$ has survived the various base changes and restrictions. Let $B$ be an abelian surface such that
$$T_\lambda(B) \ = \ \{t\in T_\lambda: \mathcal{B}^\lambda_t\cong B\}$$
has positive dimension. If $p_\lambda(\mathcal{W}_t)\subseteq  B$ does not vary with $t\in T_\lambda(B)$ there is a curve  $V\subseteq B$ such that $p_\lambda(\mathcal{W}_t)=V\subseteq B$ for all $t\in T_\lambda(B)$. Given a desingularization $\widetilde{V}$ of $V$, the abelian variety $\text{Pic}^0(\widetilde{V})$ contains an abelian subvariety isogenous to $\textup{Pic}^0(J(\mathcal{Y}_t))$ for all $t\in T_\lambda(B)$. Since a fixed abelian variety cannot contain all members of a non-isotrivial family of abelian varieties we must conclude that $p_\lambda(\mathcal{W}_t)\subseteq B$ varies with $t\in T_\lambda(B)$. This contradicts Lemma \ref{rigidhyp}.
\end{proof}

 \section{Questions and open problems}

This section is devoted to some  questions and open problems.

It is frustrating that it seems non-trivial to obtain lower bounds on $\covgon(X, Y)$ even when $X$ and $Y$ are curves. We propose
\begin{problem}
Let $X$ and $Y$ be smooth curves of genera $g_X, g_Y \ge 2$ that are Hodge-independent in the sense of \eqref{Hodge.Indpendence.Smooth.Curves}. Establish additive -- or better still, mulitiplicative  -- lower bounds on $\covgon(X, Y)$ in terms of the gonalities of $X$ and $Y$.
\end{problem}
\noi It is possible that a closer study of the bundles  $U_X$ appearing in \S 3 might be helpful for additive statements.  In a somewhat different direction, if $\gon(X)$ and $\gon(Y)$ are relatively prime and much smaller that $g_X$ and $g_Y$ can one apply the well-known argument of Castelnuovo and Severi to get a statement? John Sheridan suggests a number of other interesting questions. For example, fix a curve $Y$, an integer $d$ and a genus $g$. If $f : X \lra Y$ is a very general degree $d$ covering of $Y$ by a curve of genus $g$, does $X$ itself compute $\corrdeg(X,Y)$? Or again, can the correspondence degree between two fixed curves $X, Y$ be computed by several correspondences with essentially different numerics?

 Similarly:
\begin{problem}  Establish
better bounds (in terms of $d$ and $e$) for the order of growth of 
\[  \covgenus(X_d,Y_e) \ \ , \ \ \covgon(X_d,Y_e)\]
for very general hypersurfaces $X_d, Y_e \subseteq\PP^{n+1}$ of degrees $d$ and $e$. \end{problem}
 \noi Our sense is that  the question of establishing  multiplicative lower bounds on the joint covering gonality  has somewhat the same flavor as the problem of proving multiplicative estimates for  the covering gonality of a complete intersection of codimension $\ge 2$ in projective space.\footnote{Chen \cite{Chen} has made some striking recent progress on this question for complete intersections of codimension two and complete intersection surfaces.} In another direction, David Stapleton asks what sort of examples one can find of pairs of (special) hypersurfaces with surprisingly small correspondence degree.

 Voisin \cite{Voisin.Abel.Var} and the second author \cite{Martin} gave lower bounds on the covering gonality of very general abelian varieties of fixed dimension. It is natural to ask for an extension of those results to pairs: 
 \begin{problem}
 Let $A_1$ and $A_2$ be very general abelian varieties of dimensions $g_1$ and $g_2$. Prove an additive linear lower bound on $\covgon(A_1, A_2)$ in terms of $g_1, g_2$.
 \end{problem} 
 
In their very nice paper \cite{Chen.Stapleton}, Chen and Stapleton were able to  adapt Koll\'ar's arguments in \cite{Kollar} (involving reduction to characteristic $p > 0$) to arrive at lower bounds on measures of irrationality for complex hypersurfaces in the Fano range. This suggests
\begin{question}
Can one use these ideas to study measures of association for Fano hypersurfaces?
\end{question}

There are several  invariants of a pair of varieties $X$ and $Y$ that come to mind beyond those studied in the present paper. For example, the minimum degree of irrationality of a correspondence $Z$ between  $X$ and $Y$ measures in effect their failure to both be unirational: when $Y = \PP^n$ this degree  was studied by Yang \cite{Yang}. One might also take products with projective spaces to arrive at stable analogues of the invariants  here, and in the same spirit $\covgenus(X, Y)$ and $\covgon(X, Y)$ make perfectly good sense for pairs of varieties of different dimensions.  This suggests the somewhat vague
\begin{problem}
Explore other measures of association, and compute or estimate them in some examples. 
\end{problem}
\noi  
 Hopefully some experimentation will help to clarify what other invariants may be particularly interesting to look at. 
 
 Returning to the setting of Theorem \ref{Corr.Deg.General.X.Y.Bound}, consider a very ample line bundle $A$ on a smooth projective variety $M$ of dimension $n+1$, and let
 $X_d \in \linser{dA}$, $Y_e \in \linser{eA}$ be very general divisors. 
 \begin{question}
 Can one find conditions on $M$ and $A$ to guarantee that
 \[   
 \corrdeg(X_d, Y_e) \, \approx \, \irr(X_d) \cdot \irr(Y_e) 
 \]
 when $d, e \gg 0$?
 \end{question}
 \noi In view of  Example
 \ref{Product.CorrDeg.Example},
 it seems unlikely that this holds for every $M$ and $A$. 

 Finally, it might be interesting to explore analogues for varieties over other fields:
 \begin{problem}
 Study measures of association for varieties defined over a field $k$ other than $\CC$.
 \end{problem}
 \noi In particular, can one replace the Hodge-theoretic inputs to the results appearing above?

 \end{document}